\documentstyle{article}
\newlength{\mytopmargin}
\newlength{\myleftmargin}
\newcommand{\ka}{\kappa}
\newcommand{\al}{\alpha}
\newcommand{\de}{\delta}
\def\rlx{\relax\leavevmode}
\def\zz{\rlx\hbox{\small \sf Z\kern-.4em Z}}
\def\rr{\rlx\hbox{\scriptsize \rm I\kern-.18em R}}
\def\nn{\rlx\hbox{\rm I\kern-.18em N}}
\def\qq{\rlx\hbox{\,$\inbar\kern-.3em{\rm Q}$}}
\newcommand{\dif}[1]{\frac{\partial}{\partial #1}}
\newcommand{\bfrac}[2]{\frac{\displaystyle #1}{\displaystyle #2}}
\newcommand{\sfrac}[2]{{\textstyle\frac{#1}{#2}}}
\newcommand{\bint}{\displaystyle \int}
\newcommand{\bprod}{\displaystyle \prod}
\newcommand{\innc}[2]{\left\langle #1, #2 \right\rangle_C}
\newcommand{\innh}[2]{\left\langle #1, #2 \right\rangle_H}
\newcommand{\innl}[2]{\left\langle #1, #2 \right\rangle_L}
\newcommand{\eh}{E^{(H)}}
\newcommand{\el}{E^{(L)}}
\newcommand{\ph}{\widehat{\Phi}}

\newcommand{\wt}{\Delta(x)}

\newcommand{\et}{\widetilde{E}}
\newcommand{\shs}{\widehat{\Psi}^{\star}}
\newtheorem{lemma}{Lemma}[section]

\newtheorem{thm}[lemma]{Theorem}
\newtheorem{cor}[lemma]{Corollary}
\newtheorem{prop}[lemma]{Proposition}

\begin{document}
\noindent
\begin{center}{  \Large\bf
Non--Symmetric Jack Polynomials \\[2mm] 
and Integral Kernels }
\end{center}
\vspace{5mm}

\noindent
\begin{center} T.H.~Baker\footnote{email: tbaker@maths.mu.oz.au} and
    P.J.~Forrester\footnote{email: matpjf@maths.mu.oz.au}\\[2mm]
     {\it Department of Mathematics, University of Melbourne, \\
      Parkville, Victoria 3052, Australia}
\end{center}
\vspace{.5cm}

\begin{quote}
We investigate some properties of non-symmetric Jack, Hermite and Laguerre
polynomials which occur as the polynomial part of the eigenfunctions 
for certain Calogero-Sutherland models with exchange terms.  For the
non-symmetric Jack polynomials, the constant term normalization ${\cal
N}_\eta$ is evaluated using recurrence relations, and ${\cal N}_\eta$ is
related to the norm for the non-symmetric analogue of the power-sum inner
product. Our results for the non-symmetric Hermite and Laguerre polynomials
allow the explicit determination of the integral kernels which occur in
Dunkl's theory of integral transforms based on reflection groups of type
$A$ and $B$, and enable many analogues of properties of the classical
Fourier, Laplace and Hankel transforms to be derived. The kernels are
given as generalized hypergeometric functions based on non-symmetric
Jack polynomials. Central to our calculations is the construction of
operators $\widehat{\Phi}$ and $\widehat{\Psi}$, which act as
lowering-type operators for the non-symmetric Jack polynomials of
argument $x$ and $x^2$ respectively, and are the counterpart to the
raising-type operator $\Phi$ introduced recently by Knop and Sahi.
\end{quote}
\date{}

\setcounter{equation}{0}
\section{Introduction}

Non-symmetric Jack polynomials occur as the polynomial part of the 
eigenfunctions of the Calogero-Sutherland model on a circle with 
exchange terms. This means, in particular, that they are eigenfunctions
of the transformed Hamiltonian
\begin{eqnarray}
H^{(C)} = \sum_{j=1}^n
x_j^2 {\partial^2 \over \partial x_j^2} 
+ {2 \over \alpha}  \sum_{1 \le j < k \le n}{x_j x_k \over x_j -
x_k} \left[\left({\partial \over \partial x_j} -{\partial
\over \partial x_k} \right) - {1 - s_{jk} \over x_j - x_k} \right]
\label{h-jack}
\end{eqnarray}
Here $s_{jk}$ is the operator which acts on functions by exchanging 
the $j$'th and $k$'th coordinates. The non-symmetric Jack polynomials
$E_{\eta}(x)$, $x:=(x_1,\ldots,x_n)$ were introduced by Opdam \cite{opdam95a} 
and their properties have been expounded upon in \cite{knop96c,sahi96a}.
(Their $q$-analogues, the non-symmetric Macdonald polynomials were 
introduced in \cite{mac95} and have also received attention in the
literature \cite{cher95b,noumi96c}). In what follows, we shall mainly be
following the notation of Knop and Sahi \cite{knop96c,sahi96a}.

The $E_{\eta}(x)$ are labelled by an $n$-tuple 
$\eta=(\eta_1,\eta_2,\ldots,\eta_n)
\in\nn^n$ and are uniquely defined as being the simultaneous eigenfunctions
of the (mutually commuting) Cherednik operators $\xi_i$ defined by
\begin{equation} \label{cherednik.1}
\xi_i = \al x_i\dif{x_i} + \sum_{p<i} \frac{x_i}{x_i-x_p}(1-s_{ip})
+ \sum_{p>i} \frac{x_p}{x_i-x_p}(1-s_{ip}) +1-i
\end{equation}
and by the fact that they have an expansion of the form
\begin{equation}\label{1.3}
E_{\eta}(x) = x^{\eta} + \sum_{\nu\prec\eta}a_{\eta\nu}\,x^{\nu}\;.
\end{equation}
Here the partial order $\prec$ on $n$-tuples is defined for $\eta\neq\nu$
by 
$$
\nu\prec\eta \quad\mbox{iff}\quad \nu^+ < \eta^+ 
\quad\mbox{or in the case $\nu^+ = \eta^+$}\quad \nu < \eta
$$
where $\eta^+$ is the unique partition associated with $\eta$ obtained
from permuting its entries, and $<$ is the usual dominance order for
$n$-tuples i.e. $\nu < \eta$ iff $\sum_{i=1}^p (\eta_i - \nu_i) \geq 0$,
for all $1\leq p \leq n$. 
Indeed, we have the eigenvalue equation 
$\xi_i E_{\eta} = \bar{\eta}_i E_{\eta}$, where \cite{sahi96a}
\begin{equation}\label{e-val.1}
\bar{\eta}_i = \al\eta_i - \#\{k<i\;|\;\eta_k \geq \eta_i\} -
\#\{k>i\;|\;\eta_k > \eta_i\}
\end{equation}

The Cherednik operators are self-adjoint with respect to the inner product
\cite{cher91a}
\begin{equation}\label{inn.j}
\innc{f}{g} = {\rm C.T.}\,\left(f(x)g(x^{-1})\prod_{i\neq j}
\left(1-\frac{x_i}{x_j}\right)^{1/\al} \right)
\end{equation}
with \mbox{C.T.} meaning ``the constant term of'' in the Laurent 
polynomial expansion for $1/\al\in\nn$, and \mbox{C.T.} defined as a
Fourier integral with $x_j=e^{2\pi i\theta_j}$, $0\leq \theta_j\leq 1$,
for general $1/\al \geq 0$.
This self-adjointness, along with the form of the eigenvalues given in
(\ref{e-val.1}) implies that the non-symmetric Jack polynomials are
orthogonal with respect to the above inner product. The value of
$\innc{E_{\eta}}{E_{\eta}}$ has been computed by Opdam in \cite{opdam95a}
(and in \cite{mac95,cher95b} for the $q$-case). One of the new results of this
work, given in Section 2, is a different evaluation of
$\innc{E_{\eta}}{E_{\eta}}$,
utilizing the operator $\Phi$ introduced in \cite{knop96c}.

Recently we began a study of the
eigenfunctions of the Calogero-Sutherland model in an external harmonic
potential with exchange terms, associated with the roots systems of
type $A$ and $B$ \cite{forr96c}. The transformed Hamiltonians take 
the form 
\begin{eqnarray}
{H}^{(H)} & := & \sum_{j=1}^n
\left( {\partial^2 \over \partial x_j^2 } - 2 x_j {\partial
\over \partial x_j } \right) + {2 \over \alpha} \sum_{j < k}
{1 \over x_j - x_k} \left[ \left( {\partial \over \partial x_j } -
{\partial \over \partial x_k } \right) - {1 - s_{jk} \over x_j - x_k}
\right] \label{h-herm} \\
{H}^{(L)} & := & \sum_{j=1}^n \left(
{\partial^2 \over \partial x_j^2 } + \left(\frac{2a+1}{x_j} - 2x_j\right)
{\partial \over \partial x_j } \right)\nonumber \\
&& + {4 \over \alpha} \sum_{j < k} {1 \over x^2_j - x^2_k}
\left[ \left( x_j {\partial \over \partial x_j } -
x_k  {\partial \over \partial x_k } \right) -
{x^2_j + x^2_k \over x^2_j - x^2_k} (1 -s_{jk}) \right],
\label{h-lag}
\end{eqnarray}
where to obtain (\ref{h-lag}) we have set $y_j=x_j^2$ in
\cite[(1.17)]{forr96c} and multiplied through by 4,
and their eigenfunctions are called non-symmetric Hermite 
and Laguerre polynomials, denoted $\eh_\eta(x)$ and $\el_\eta(x^2)$
respectively. In ref.~\cite{forr96c} we constructed a set
of commuting operators for $H^{(H)}$, which have $\{E_\eta^{(H)}\}$
as simultaneous eigenfunctions, and which are self adjoint with respect
to the inner product
\begin{equation}
\innh{f}{g} := \prod_{l=1}^n \int_{-\infty}^\infty
dx_l \, e^{-x_l^2} \prod_{1 \le j < k \le n} |x_k - x_j|^{2/\alpha} 
\,f\, g
\label{inn.h} 
\end{equation}
A consequence of this construction is the orthogonality relation  
$\langle E_\eta^{(H)}, E_\nu^{(H)} \rangle_{H} = {\cal N}_\eta^{(H)}
\delta_{\eta,\nu}$,
where ${\cal N}_\eta^{(H)}$ is the norm.
Similarly, we constructed a set of commuting operators
for $H^{(L)}$, which are self adjoint with respect to
\begin{equation} \label{inn.l}
\innl{f}{g} := 2^n \prod_{l=1}^n \int_{-\infty}^\infty
dx_l \,
e^{-x_l^2} |x_l|^{2a+1} \prod_{1 \le j < k \le n} |x_k^2 - x_j^2|^{2/\alpha}
 f(x^2) g(x^2)
 \end{equation}

In this paper, we continue our study of the non-symmetric Jack, Hermite
and Laguerre polynomials. In Section 2
we review some results concerning the non-symmetric Jack polynomials,
in particular those due to Knop and Sahi \cite{knop96c,sahi96a}, which
are relevant to the calculations in later sections. Using these results, we
provide a new proof of the evaluation of the norm $\innc{E_{\eta}}
{E_{\eta}}$. We also give a non-symmetric analogue of a generalization
due to Kadell \cite{kad93a} of the Morris constant term identity 
\cite{morrthesis}. In the course of deriving this result we obtain 
a formula relating the norm $\innc{E_{\eta}}{E_{\eta}}$ and the norm 
$\langle E_{\eta}\,,\,E_{\eta} \rangle$, where 
$\langle \cdot\,,\,\cdot\rangle$ is the non-symmetric analogue of
the power-sum inner product specified in \cite{dunkl96a,sahi96a}.

In Section 3, we turn our attention to the non-symmetric Hermite
polynomials $\eh_{\eta}$, beginning with a brief review of the
pertinent results in \cite{forr96c}. We then
proceed to construct certain operators which are sufficient to generate
all $\eh_{\eta}$ by recurrence. This enables us to compute the norm
${\cal N}_{\eta}^{(H)}$ for non-symmetric Hermite polynomials.
Our attention is then directed towards a construction of Dunkl's 
\cite{dunkl91a,dunkl92a} 
non-symmetric kernel ${\cal K}_A(x;y)$ (an analogue of the (symmetric)
generalized hypergeometric series ${}_0{\cal F}_0(x;y)$ \cite{yan92b}),
which is used to define generalizations of the Fourier and Laplace
transforms.
This kernel also allows us to derive an exponential formula, a generating
function, and integral formulae for the non-symmetric Hermite polynomials
in complete analogy with the symmetric case \cite{forr96a}. Indeed,
we show that ${}_0{\cal F}_0(x;y)$ can be constructed from ${\cal K}_A(x;y)$
by symmetrization. The analysis of Section 3
is repeated, albeit more succinctly, in Section 4 for the Laguerre case.

\setcounter{equation}{0}
\section{The Jack case}

We begin by reviewing some of the results in \cite{knop96c,sahi96a}. A
fundamental result concerns the action of the elementary transpositions
$s_i:=s_{i,i+1}$ on the non-symmetric Jack polynomials, which is given
by
\begin{lemma}\label{sahi.1}
Let $\bar{\eta}_i$ be the eigenvalue of $\xi_i$ on the non-symmetric Jack
polynomial $E_{\eta}$, and let $\de_{i,\eta}:=\bar{\eta}_i -
\bar{\eta}_{i+1}$.
Then the action of $s_i$ is given by
$$
s_i\,E_{\eta} = \left\{ \begin{array}{ll}
\bfrac{1}{\de_{i,\eta}}\,E_{\eta} + \left( 1- \bfrac{1}{\de_{i,\eta}^2}
\right) E_{s_i\eta} & \eta_i > \eta_{i+1} \\
E_{\eta} & \eta_i = \eta_{i+1} \\
\bfrac{1}{\de_{i,\eta}}\,E_{\eta} + E_{s_i\eta} & \eta_i < \eta_{i+1}
\end{array} \right.
$$
\end{lemma}
This is a consequence of the following relations between the Cherednik
operators and the transpositions $s_i$
\begin{equation}\label{hecke}
\xi_i s_i - s_i \xi_{i+1} = 1, \qquad \xi_{i+1}s_i -s_i \xi_i = -1, \qquad
{}[\xi_i,s_j] = 0, \quad j\neq i, i+1
\end{equation}

Knop and Sahi also introduced a remarkable operator
$\Phi$, defined by
\begin{equation}\label{phi-op}
\Phi = x_n\, s_{n-1}\cdots s_2\, s_1 = s_{n-1}\cdots s_i \:x_i \:
s_{i-1}\cdots s_1
\end{equation}
which enjoys the following properties
\begin{lemma}\label{sahi.2}
\begin{eqnarray*}
\xi_j\,\Phi &=& \Phi\,\xi_{j+1} \hspace{2cm} 1\leq j \leq n-1 \\
\xi_n\,\Phi &=& \Phi\, (\xi_1 +\al) \\
\Phi\,E_{\eta}(x) &=& E_{\Phi\eta} (x)
\end{eqnarray*}
where $\Phi\eta := (\eta_2,\eta_3,\ldots,\eta_n,\eta_1+1)$
\end{lemma}
As noted in \cite{knop96c}, these results imply that the operators
$s_i$, $1\leq i\leq n$ and $\Phi$ are sufficient to generate all
$E_{\eta}$. As an application, Sahi \cite{sahi96a} was subsequently able
to evaluate $E_{\eta}$ at the point $x_1=x_2=\cdots=x_n=1$. To write
down this result, which is required below, we need some additional 
notation. For a node $s=(i,j)$ in an $n$-tuple
$\eta$, define the arm length $a(s)$, arm colength $a'(s)$, leg length
$l(s)$ and  leg colength $l'(s)$ by
\begin{eqnarray}
a(s)= \eta_i - j && l(s) = \#\{k>i|j\leq \eta_k\leq\eta_i\} \;+\;
\#\{k<i|j\leq \eta_k+1\leq\eta_i\} \nonumber\\
a'(s)=j - 1 && l'(s) = \#\{k>i| \eta_k > \eta_i\} \;+\;
\#\{k<i|\eta_k\geq\eta_i\}  \label{guion}
\end{eqnarray}
Using these, define constants $d_{\eta}:=\prod_{s\in\eta} d(s)$,
$d'_{\eta}:=\prod_{s\in\eta} d'(s)$ and $e_{\eta}:=\prod_{s\in\eta} e(s)$
where
$$
d'(s) := \al(a(s)+1) + l(s)\hspace{2cm} e(s):= \al(a'(s)+1) + n-l'(s)
$$
and $d(s):=d'(s)+1$. These have the following important properties
\begin{lemma}\label{recur.1}
We have
\begin{eqnarray*}
\frac{d_{\Phi\eta}}{d_{\eta}} &=& \frac{e_{\Phi\eta}}{e_{\eta}}
= \bar{\eta}_1 +\al +n \hspace{2cm} \frac{d'_{\Phi\eta}}{d'_{\eta}} = 
\bar{\eta}_1 +\al +n-1 \qquad\mbox{for all\  $\eta$}\\
e_{s_i\eta} &=& e_{\eta} \qquad \frac{d_{s_i\eta}}{d_{\eta}} = 
\frac{\de_{i,\eta}+1}{\de_{i,\eta}} \qquad 
\frac{d'_{s_i\eta}}{d'_{\eta}} = \frac{\de_{i,\eta}}{\de_{i,\eta}-1}
\qquad \mbox{for $\eta_i > \eta_{i+1}$}
\end{eqnarray*}
A similar relation follows in the case $\eta_i<\eta_{i+1}$ after
noting that $\de_{i,s_i\eta} = - \de_{i,\eta}$.
\end{lemma}
Using Lemmas \ref{sahi.1}, \ref{sahi.2} and \ref{recur.1}, Sahi showed that
\begin{equation}
E_{\eta}(1^n) = e_{\eta}/d_{\eta}
\end{equation}
by showing that both sides of this equation satisfy the same recursions
via the operators $s_i$ and $\Phi$ (see \cite{cher95b} for another proof
of this result).

\subsection{Calculation of $\langle E_\eta, E_\eta \rangle_C$}

Let us now show that a similar idea works for the calculation of the
norm of the non-symmetric Jack polynomials with respect to the inner
product (\ref{inn.h}).

\begin{prop}
In the case where $k=1/\al\in\zz^+$
\begin{equation}\label{socavar}
\innc{E_{\eta}}{E_{\eta}} = \prod_{1\leq i < j \leq n} \prod_{p=0}^{k-1}
\left( \frac{k(\bar{\eta}_j - \bar{\eta}_i) + p}
{k(\bar{\eta}_j - \bar{\eta}_i) - p -1} \right)^{\epsilon(
\bar{\eta}_j - \bar{\eta}_i)}
\end{equation}
where $\epsilon(x) = 1$ for $x>0$, and $\epsilon(x) = -1$ for $x\leq 0$
\end{prop}
{\it Proof.}\quad First, note 
that the transpositions $s_i$ are hermitian w.r.t.~the inner product 
defined by (\ref{inn.j}). Thus from the definition (\ref{phi-op})
$\innc{\Phi\,f}{\Phi\,g} = \innc{f}{g}$, 
and so $\Phi$ is an isometry.  Thus from Lemma \ref{sahi.2} we have
\begin{equation}\label{hito.1}
\innc{E_{\Phi\eta}}{E_{\Phi\eta}} = \innc{\Phi\,E_{\eta}}{\Phi\,E_{\eta}}
= \innc{E_{\eta}}{E_{\eta}}
\end{equation}
Also, for the case $\eta_i < \eta_{i+1}$, from Lemma \ref{sahi.1} 
we have $E_{s_i\eta} =
s_i\,E_{\eta} - \de_{i,\eta}^{-1} E_{\eta}$ so that
\begin{eqnarray}
\innc{E_{s_i\eta}}{E_{s_i\eta}} &=&  \innc{(s_i -\de_{i,\eta}^{-1})E_{\eta}}
{(s_i -\de_{i,\eta}^{-1})E_{\eta}}\nonumber\\
&=& (1+\de_{i,\eta}^{-2})\innc{E_\eta}{E_\eta} -2\de_{i,\eta}^{-1}
\innc{E_{\eta}}{s_i\,E_{\eta}} \nonumber\\
&=& (1-\de_{i,\eta}^{-2})\:\innc{E_\eta}{E_\eta} \label{hito.2}
\end{eqnarray}
where in obtaining the last line we have used Lemma \ref{sahi.1} in the
second term of the previous line, and used the fact that for $\eta_i
\neq \eta_{i+1}$, $E_{\eta}$ and $E_{s_i\eta}$ are orthogonal. Equation
(\ref{hito.2}) immediately implies an equivalent result in the case 
$\eta_i > \eta_{i+1}$, namely
$$
\innc{E_{s_i\eta}}{E_{s_i\eta}} =
\left(1-\de_{i,\eta}^{-2}\right)^{-1}\:\innc{E_\eta}{E_\eta}
$$
through the obvious change of variables $\eta\rightarrow s_i\eta$
(recall $\de_{i,s_i\eta} = - \de_{i,\eta}$). It thus remains to show
that the right hand side of (\ref{socavar}), ${\rm RHS}_{\eta}$ say, 
obeys the same recursion
relations (\ref{hito.1}) and (\ref{hito.2}) and that both sides have
the same evaluation in the trivial case $\eta=0$. For the latter property
note that then $\bar{\eta}_i = 1-i$ and so (\ref{socavar}) reduces to
$$
\innc{1}{1} = \prod_{1\leq i<j\leq n}\prod_{p=0}^{k-1}
\left(\frac{k(i-j) -p-1}{k(i-j)+p} \right)
$$
which is a well-known constant term identity (see for example \cite{mac}).

Turning to (\ref{hito.1}), first note from the definition (\ref{e-val.1})
that if $\zeta=\Phi\eta$, then 
\begin{equation}\label{necio}
(\bar{\zeta}_1,\bar{\zeta}_2,\ldots,\bar{\zeta}_n) =
(\bar{\eta}_2,\bar{\eta}_3,\ldots,\bar{\eta}_n,\bar{\eta}_1+\al)
\end{equation}
Thus we have
\begin{eqnarray}\lefteqn{ \mbox{RHS}_{\eta}
=\prod_{1\leq i < j \leq n} \prod_{p=0}^{k-1}
\left( \frac{k(\bar{\zeta}_j - \bar{\zeta}_i) + p}
{k(\bar{\zeta}_j - \bar{\zeta}_i) - p -1} \right)^{\epsilon(
\bar{\zeta}_j - \bar{\zeta}_i)} =} && \nonumber\\
&&\prod_{2\leq i < j \leq n} \prod_{p=0}^{k-1}
\left( \frac{k(\bar{\eta}_j - \bar{\eta}_i) + p}
{k(\bar{\eta}_j - \bar{\eta}_i) - p -1} \right)^{\epsilon(
\bar{\eta}_j - \bar{\eta}_i)}\;
\prod_{i=2}^n\prod_{p=0}^{k-1}
\left( \frac{k(\bar{\eta}_1+\al - \bar{\eta}_i) + p}
{k(\bar{\eta}_1+\al - \bar{\eta}_i) - p -1} \right)^{\epsilon(
\bar{\eta}_1+\al - \bar{\eta}_i)} \label{osado}
\end{eqnarray}
If we can show that $\epsilon(\bar{\eta}_1+\al - \bar{\eta}_i) = 
- \epsilon(\bar{\eta}_i - \bar{\eta}_1)$ for all $2\leq i\leq n$,
then the second double product on the right hand side of (\ref{osado})
can be rewritten as
$$
\prod_{i=2}^n\prod_{p=0}^{k-1}
\left( \frac{-[k(\bar{\eta}_i - \bar{\eta}_1) - p -1]}
{-[k(\bar{\eta}_i - \bar{\eta}_1) + p]} \right)^{-\epsilon(
\bar{\eta}_i - \bar{\eta}_1)}
= \prod_{i=2}^n\prod_{p=0}^{k-1}
\left( \frac{k(\bar{\eta}_i - \bar{\eta}_1) + p }
{k(\bar{\eta}_i - \bar{\eta}_1) - p -1} \right)^{\epsilon(
\bar{\eta}_i - \bar{\eta}_1)}
$$
which when reinserted back into (\ref{osado}) yields the required  
equality. 

To show that $\epsilon(\bar{\eta}_1+\al - \bar{\eta}_i) = 
- \epsilon(\bar{\eta}_i - \bar{\eta}_1)$, we consider two cases. In
the first case, suppose $\epsilon(\bar{\eta}_1+\al - \bar{\eta}_i) =
-1$, i.e. $\bar{\eta}_1+\al \leq \bar{\eta}_i$. Then $\bar{\eta}_i -
\bar{\eta}_1\geq \al >0$, so that $\epsilon(\bar{\eta}_i - \bar{\eta}_1)
=1$ as required. In the second case, where $\epsilon(\bar{\eta}_1+\al - 
\bar{\eta}_i) =1$, we have $\bar{\eta}_1+\al > \bar{\eta}_i$. Note
that we can always write $\bar{\eta}_1 - \bar{\eta}_i = a\al+b$ for some
integers $a$ and $b$. Then using the fact that $\al=1/k$, we have
$\bar{\eta}_1+\al > \bar{\eta}_i \Leftrightarrow bk>-a-1$. But $b$, $k$
and $a$ are all integers, so that this later inequality is equivalent to 
$bk\geq -a$ which in turn is equivalent to $\bar{\eta}_i - 
\bar{\eta}_1 \leq 0$, and hence $\epsilon(\bar{\eta}_i - \bar{\eta}_1) =
-1$.

Turning to (\ref{hito.2}), note that if $\nu =s_i\eta$, then 
$$
(\bar{\nu}_1,\ldots,\bar{\nu}_{i-1},
\bar{\nu}_{i},\bar{\nu}_{i+1},\bar{\nu}_{i+2},\ldots,\bar{\nu}_{n}) =
(\bar{\eta}_1,\ldots,\bar{\eta}_{i-1},
\bar{\eta}_{i+1},\bar{\eta}_{i},\bar{\eta}_{i+2},\ldots,
\bar{\eta}_{n})
$$
Moreover, when $\eta_i < \eta_{i+1}$, then $\bar{\eta}_{i+1} - \bar{\eta}_i
\geq \al(\eta_{i+1}-\eta_i)+1 > 0$ with a similar result occuring when
$\eta_i > \eta_{i+1}$, so that in all cases $\bar{\eta}_{i+1} -
\bar{\eta}_i \neq 0$. The upshot of all of this is that we always have
$\epsilon(\bar{\eta}_i-\bar{\eta}_{i+1}) = - 
\epsilon(\bar{\eta}_{i+1}-\bar{\eta}_i)$. From these considerations we
see that
\begin{eqnarray*} 
\prod_{p=0}^{k-1} \left(\frac{k(\bar{\nu}_{i+1} - \bar{\nu}_{i})+p}
{k(\bar{\nu}_{i+1} - \bar{\nu}_{i})-p-1}\right)^{\epsilon(
\bar{\nu}_{i+1} - \bar{\nu}_{i}) } 
= \prod_{p=0}^{k-1} \left(\frac{k(\bar{\eta}_{i+1} - \bar{\eta}_{i})+p+1}
{k(\bar{\eta}_{i+1} - \bar{\eta}_{i})-p}\right)^{\epsilon(
\bar{\eta}_{i+1} - \bar{\eta}_{i}) }  \hspace{2cm}\\
=\left(1-\de_{i,\eta}^{-2} \right)^{\epsilon(
\bar{\eta}_{i+1} - \bar{\eta}_{i}) }\;
\prod_{p=0}^{k-1} \left(\frac{k(\bar{\eta}_{i+1} - \bar{\eta}_{i})+p}
{k(\bar{\eta}_{i+1} - \bar{\eta}_{i})-p-1}\right)^{\epsilon(
\bar{\eta}_{i+1} - \bar{\eta}_{i}) }
\end{eqnarray*}
which in turn implies that the right hand side of (\ref{socavar}) 
satisfies the relation (\ref{hito.2}). \hfill $\Box$

\subsection{Generalization of the Macdonald-Kadell-Kaneko
integral}
The Macdonald-Kadell-Kaneko integral, first conjectured by Macdonald
\cite{mac87} and subsequently proved by Kadell \cite{kad93a} and
Kaneko \cite{kaneko93a}, relates the Selberg integral \cite{selberg}
to the Jack polynomial. Here we will derive the analogue of this
result, relating the Selberg integral to the non-symmetric Jack 
polynomial.

Our derivation  is based on the generalized binomial theorem
\cite{stan89a,kad93a,kaneko93a} 
\begin{equation}\label{rebanada}
\prod_{i=1}^n \frac{1}{(1-x_i)^r} = \sum_{\kappa} \frac{\al^{|\ka|}
[r]^{(\al)}_{\ka}}{j_{\kappa}} \; J^{(\al)}_{\kappa}(x)
\end{equation}
where 
\begin{equation} \label{blitz}
[r]^{(\al)}_{\ka} := \prod_{j=1}^n \frac{\Gamma(r-(j-1)/\al +\ka_j)}
{\Gamma(r-(j-1)/\al )}
\end{equation}
Let us introduce the non-symmetric Jack polynomials with a different
normalization namely $F_{\eta}:=d_{\eta}E_{\eta}$, and define constants
$f_\eta:=d_{\eta}d'_{\eta}$. 
{}From \cite{sahi96a,mac95} we know that
\begin{equation} \label{abolladura}
\frac{1}{j_\ka} J^{(\al)}_{\kappa}(x) = \sum_{\eta : \eta^+ = \kappa} 
\frac{1}{f_{\eta}} F_{\eta}(x)
\end{equation}
Hence from (\ref{rebanada}) we have
\begin{equation} \label{lonja}
\prod_{i=1}^n \frac{1}{(1-x_i)^r} = \sum_{\eta} \frac{\al^{|\eta|}
[r]^{(\al)}_{\eta^+}}{f_{\eta}} \; F_{\eta}(x)
\end{equation}
Recall that $\{F_{\eta}\}$ is an orthogonal set of functions with respect
to the constant-term inner product (\ref{inn.j}), and that it constitutes
a basis for analytic functions. Hence letting $\wt:=\prod_{j\neq k}
(1-\frac{x_j}{x_k})^{2/\al}$ we can write
\begin{equation} \label{rodaja}
\prod_{i=1}^n \frac{1}{(1-x^{-1}_i)^r} = \sum_{\eta} \frac
{{\rm C.T.}(\prod_{i=1}^n (1-x^{-1}_i)^{-r} F_{\eta}(x) \wt ) }
{{\rm C.T.}(F_{\eta}(x^{-1}) F_{\eta}(x) \wt) }  \; F_{\eta}(x^{-1})
\end{equation}
Comparing (\ref{lonja}) and (\ref{rodaja}) gives, after replacing $1/x$
by $x$ in (\ref{rodaja}),
\begin{equation}\label{zapato}
{\rm C.T.}(\prod_{i=1}^n (1-x^{-1}_i)^{-r} F_{\eta}(x)\wt) =
\frac{\al^{|\eta|} [r]^{(\al)}_{\eta^+}}{f_{\eta}}
{\rm C.T.}(F_{\eta}(x^{-1})F_{\eta}(x)\wt)
\end{equation}

Our next task is to manipulate the left hand side of (\ref{zapato}) so
that $(1-x_i^{-1})^{-r}$ is replaced by $(1-x_i)^a(1-x_i^{-1})^b$. We
require

\begin{lemma} \label{yep}
We have
$$
x^p\:E_{\eta}(x) = E_{\eta+p}(x)
$$
where $\eta+p:=(\eta_1+p,\eta_2+p,\ldots,\eta_n+p)$, and
$x^p:=(x_1x_2\cdots x_n)^p$.
\end{lemma}
{\it Proof.}\quad  Using the Cherednik operators (\ref{cherednik.1})
and the corresponding eigenvalue equation for $E_{\eta}(x)$, we have
$$
\xi_j\; (x^p\,E_{\eta}) = (\bar{\eta}_j +\al p)
x^p\,E_{\eta} .
$$
However, from the definition (\ref{e-val.1}) we have
$\bar{\eta}_j +\al p = \overline{(\eta+p)}_j$ so that $x^p\,E_{\eta}$
must be a constant multiple of $E_{\eta+p}$. Examination of the
leading terms shows that this constant is $1$ and the result then
follows. \hfill $\Box$

\vspace{.2cm}

Consider (\ref{zapato}) with $\eta$ replaced by $\eta +a$. Now
$$
F_{\eta+a}(x) = d_{\eta+a} E_{\eta+a}(x) = d_{\eta+a}x^a E_{\eta}(x)
$$
by Lemma \ref{yep}. Also, set $r=-a-b$ and note that
$$
x_i^a(1-\frac{1}{x_i})^{a+b} = (-1)^a(1-x_i)^a(1-\frac{1}{x_i})^b
$$
This gives
\begin{equation}
{\rm C.T.}\left(\prod_{i=1}^n(1-x_i)^a(1-\frac{1}{x_i})^b E_{\eta}(x)
\wt \right)  
= (-\al)^{an} \al^{|\eta|}[-a-b]^{(\al)}_{\eta^++a}
\frac{d_{\eta+a}}{f_{\eta+a}} \;{\rm C.T.}\left(E_{\eta}(x^{-1})
E_{\eta}(x) \wt \right) \label{acicalarse}
\end{equation}
The dependence on $a$ in $d_{\eta+a}/f_{\eta+a}$ can be determined by
using

\begin{lemma} \label{reverse}
We have
$$
E_{\eta}\left(\frac{1}{x}\right) = E_{-\underline{\eta}}(\underline{x})
$$
where $\underline{x}:=(x_n,x_{n-1},\ldots,x_1)$, 
$\underline{\eta}:=(\eta_n,\eta_{n-1},\ldots,\eta_1)$ and $E_{-
\underline{\eta}}$
is interpreted according to Lemma \ref{yep}.
\end{lemma}
{\it Proof.}\quad Let $y=1/x$. Then 
$$
\xi^{(y)}_i\;E_{\eta}(y) = \bar{\eta}_i \;E_{\eta}(y) .
$$
Since $\dif{y_i} = -x_i^2 \dif{x_i}$, we have
\begin{eqnarray*}
\xi^{(y)}_i &=& -\al x_i\dif{x_i} -\sum_{p<i} \frac{x_p}{x_i-x_p}(1-s_{ip})
- \sum_{p>i} \frac{x_i}{x_i-x_p}(1-s_{ip}) + 1-i \\
&=& -\al \hat{D}_i^{(x)}
\end{eqnarray*}
where $\hat{D}_i$ is Hikami's version \cite{hik96a} of the Cherednik
operator (see ref.~\cite[eq. (2.2) and the remarks which follow]{forr96c}).
Hence
$$
\al \hat{D}_i^{(x)}\;E_{\eta}\left(\frac{1}{x}\right) = -\bar{\eta}_i\;
E_{\eta}\left(\frac{1}{x}\right) .
$$
But the leading term of $E_{\eta}\left(\frac{1}{x}\right)$ is $(\frac{1}
{x})^{\eta}$ as is the leading term of
$E_{-\underline{\eta}}(\underline{x})$
which is also an eigenfunction of $\hat{D}_i^{(x)}$. Hence the
equality. \hfill $\Box$

\vspace{2mm}\noindent
Now
\begin{eqnarray*} 
{\rm C.T.}\left(\prod_{i=1}^n(1-x_i)^a(1-\frac{1}{x_i})^b E_{\eta}
(x) \wt\right) 
&=& {\rm C.T.}\left(\prod_{i=1}^n(1-\frac{1}{x_i})^a(1-x_i)^b 
E_{\eta}(\underline{x}^{-1})\wt \right)  \\
&=& {\rm C.T.}\left(\prod_{i=1}^n(1-\frac{1}{x_i})^a(1-x_i)^b 
E_{-\underline{\eta}}(x) \wt\right)  
\end{eqnarray*}
where to obtain the second line we have used the invariance of
the constant-term operation under $x_i\rightarrow 1/x_i$ and
permutations of the variables, along with the symmetry of the terms
appearing in the integrand (excluding of course $E_{\eta}$), and
to get the last line we have used Lemma \ref{reverse}. Hence
(\ref{acicalarse}) is unchanged if we interchange $a$ and $b$ and
replace $\eta$ by $-\underline{\eta}$. Equating the corresponding right hand
sides of (\ref{acicalarse}) gives
\begin{equation} \label{inti}
(-\al)^{an}  \al^{|\eta|}[-a-b]^{(\al)}_{\eta^++a}
\frac{d_{\eta+a}}{f_{\eta+a}} =
(-\al)^{bn}  \al^{-|\eta|}[-a-b]^{(\al)}_{-\underline{\eta}^++b}
\frac{d_{-\underline{\eta}+b}}{f_{-\underline{\eta}+b}}
\end{equation}
where we have used the fact that 
$$
{\rm C.T.}\left(E_{\eta}(x^{-1})E_{\eta}(x)\wt\right) =
{\rm C.T.}\left(E_{-\underline{\eta}}(x^{-1})E_{-\underline{\eta}}(x)\wt\right).
$$
Setting $a=0$ in (\ref{inti}) gives 
\begin{equation} \label{inti.2}
(-\al)^{bn}  \al^{-|\eta|} \frac{d_{-\underline{\eta}+b}}
{f_{-\underline{\eta}+b}}
=\al^{|\eta|}\frac{[-b]^{(\al)}_{\eta^+}}{[-b]^{(\al)}_{-\underline{\eta}^++b}}
\frac{d_{\eta}}{f_{\eta}}
\end{equation}
Substituting (\ref{inti.2}) into (\ref{acicalarse}) with $a$ and $b$
interchanged on the right hand side and $\eta$ replaced by $-\underline{\eta}$
gives
\begin{equation}
{\rm C.T.}\left(\prod_{i=1}^n(1-x_i)^a(1-\frac{1}{x_i})^b
E_{\eta}(x) \wt\right) 
= \al^{|\eta|}\frac{[-b]^{(\al)}_{{\eta}^+}
[-a-b]^{(\al)}_{-\underline{\eta}^++b}}{[-b]^{(\al)}_{-\underline{\eta}^++b}}
\frac{d_{\eta}}{f_{\eta}}\; N_{\eta} \label{hoja}
\end{equation}
where ${\cal N}_{\eta}:={\rm C.T.}\left(E_{\eta}(x^{-1})
E_{\eta}(x)\wt\right)$. In (\ref{hoja}), the dependence on $a$ and $b$ 
on the right hand side is explicitly displayed by virtue of
(\ref{blitz}). This allows the $a$, $b\rightarrow\infty$ asymptotics
to be computed, which leads to a relationship between $E_{\eta}(1^n)$
and $d_{\eta}\;{\cal N}_{\eta}/f_{\eta}$ and also allows a simplification 
of the right hand side of (\ref{hoja}). It is convenient to first 
take the ratio of (\ref{hoja}) to that obtained upon setting
$\eta=0$ :
\begin{eqnarray}
\frac{{\rm C.T.}\left(\prod_{i=1}^n(1-x_i)^a(1-\frac{1}{x_i})^b
E_{\eta}(x_1,\ldots,x_n) \wt\right) }
{{\rm C.T.}\left(\prod_{i=1}^n(1-x_i)^a(1-\frac{1}{x_i})^b
\wt\right) } \hspace{4cm}\nonumber\\
= \al^{|\eta|}\frac{[-b]^{(\al)}_{\eta^+}
[-a-b]^{(\al)}_{-\underline{\eta}^++b}}{[-b]^{(\al)}_{-\underline{\eta}^++b}}
\frac{d_{\eta}}{f_{\eta}} \;\frac{[-b]^{(\al)}_b}{[-a-b]^{(\al)}_b}\;
\frac{{\cal N}_{\eta}}{{\cal N}_0} \label{ash}
\end{eqnarray}
Now set $a=b$ and take the limit $a\rightarrow\infty$ in
(\ref{ash}). On the left hand side the maximum contribution comes from
the neighbourhood of $x_1=\cdots x_n =-1$, and so the ratio is equal to
$E_{\eta}(-1^n)$ while on the right hand side, since
\begin{eqnarray*}
{}[-r]^{(\al)}_{\ka} &=& \prod_{j=1}^n \frac{\Gamma(-r-(j-1)/\al
+\ka_j)}{\Gamma(-r-(j-1)/\al )} \\
&=& \prod_{j=1}^n \frac{\sin\pi(-r-(j-1)/\al)\;\Gamma(1+r+(j-1)/\al)}
{\sin\pi(-r-(j-1)/\al+\ka_j)\;\Gamma(1+r+(j-1)/\al -\ka_j)} \\
&=& (-1)^{|\ka|} \prod_{j=1}^n \frac{\Gamma(1+r+(j-1)/\al)}
{\Gamma(1+r+(j-1)/\al -\ka_j)} 
\end{eqnarray*}
and $\Gamma(x+u)/\Gamma(x) \sim x^u$ as $x\rightarrow\infty$, we see
that
\begin{eqnarray*}
[-b]^{(\al)}_{{\eta}^+} &\sim & (-1)^{|\eta|}\,b^{|\eta|} \\
\frac{[-a-b]^{(\al)}_{-\underline{\eta}^++b}}
{[-a-b]^{(\al)}_{b}} &=& (-1)^{|\eta|}\; \frac{1}{[1+a+(n-1)/\al]^{(\al)}_
{\eta^+}} \sim (-1)^{|\eta|}\,a^{-|\eta|} \\
\frac{[-b]^{(\al)}_{b}}{[-b]^{(\al)}_{-\underline{\eta}^++b}} &=&
(-1)^{|\eta|}\;\prod_{j=1}^n \frac{\Gamma(1+(j-1)/\al)}
{\Gamma(1+(j-1)/\al + \eta^+_{n+1-j})}
\end{eqnarray*}
Hence with $a=b$, in the limit $a\rightarrow\infty$ (\ref{ash})
gives
$$
E_{\eta}(-1^n) = (-1)^{|\eta|} \left(
\prod_{j=1}^n \frac{\Gamma(1+(j-1)/\al)}
{\Gamma(1+(j-1)/\al +\eta^+_{n+1-j})} \right) \al^{|\eta|}\frac{
d_{\eta}}{e_{\eta}}\;\frac{{\cal N}_{\eta}}{{\cal N}_0}
$$
and so
\begin{equation} \label{barbilla}
\al^{|\eta|}\frac{d_{\eta}}{f_{\eta}}\;\frac{{\cal N}_{\eta}}{{\cal N}_0}
= E_{\eta}(1^n)\left(\prod_{j=1}^n \frac{\Gamma(1+(j-1)/\al)}
{\Gamma(1+(j-1)/\al +\eta^+_{n+1-j})} \right)^{-1} .
\end{equation}
This relates $d_{\eta}\;{\cal N}_{\eta}/f_{\eta}$ to $E_{\eta}(1^n)$ 
and also relates the constant term norm ${\cal N}_{\eta}$ to the norm
$\langle E_{\eta}\,,\,E_{\eta}\rangle$, where 
$\langle \cdot\,,\,\cdot\rangle$ is the non-symmetric analogue of the
power-sum inner product specified in \cite{dunkl96a,sahi96a} (for the
non-symmetric Macdonald polynomials, the result analogous to 
(\ref{barbilla}) has been given recently by Mimachi and Noumi
\cite{noumi96c}).

We can therefore rewrite (\ref{ash}) as
\begin{equation}
\frac{{\rm C.T.}\left(\prod_{i=1}^n(1-x_i)^a(1-\frac{1}{x_i})^b
E_{\eta}(x_1,\ldots,x_n) \wt\right) }
{{\rm C.T.}\left(\prod_{i=1}^n(1-x_i)^a(1-\frac{1}{x_i})^b
\wt\right) } 
=E_{\eta}(1^n)\;\frac{[-b]^{(\al)}_{\eta^+}}{[1+a+(n-1)/
\al]^{(\al)}_{\eta^+}} \label{sgarden}
\end{equation}
Note that by multiplying both sides by $d_{\eta}/e_{\eta}$ and 
summing over the distinct permutations of $\ka=\eta^+$ using
(\ref{abolladura}), we get back 
(\ref{sgarden}) with $E_{\eta}$ replaced by $J_{\ka}^{(\al)}$ on both
sides, which is the formula of Kadell \cite{kad93a}.

To obtain the integration formula of Macdonald, Kadell and Kaneko, we
make use of the following lemma, proved in \cite{forr95b}

\begin{lemma}\label{plem}
For Re$(\epsilon)$ large enough so that the right hand side exists,
$$
\left(\frac{\pi}{\sin\pi\epsilon}\right)^n\;\prod_{l=1}^n
\int^{1/2}_{-1/2} d\theta_l\;e^{2\pi i\theta_l\epsilon}\;
f\left(-e^{2\pi i\theta_1},\ldots,-e^{2\pi i\theta_n} \right)
= \prod_{l=1}^n \int^{1}_0 dt_l\; t_l^{-1+\epsilon}\;
f(t_1,\ldots,t_n)
$$
where $f$ is a Laurent polynomial.
\end{lemma}

\vspace{2mm}\noindent
Notice from the derivation of (\ref{sgarden}) that $a+b\:(\:=-r)$ is
arbitrary, as is $b$. From the symmetry relation with respect
to interchanging $a$ and $b$, it follows that $a$ is arbitrary as well.
Now
$$
\Delta(x) = \left(\prod_{1\leq j<k\leq n}\left(1-\frac{x_k}{x_j}\right)
\left(1-\frac{x_j}{x_k}\right) \right)^{1/\al}\!\!\! =
(-1)^{n(n-1)/2\al} \prod_{j=1}^n x_j^{-(n-1)/\al}\!
\prod_{1\leq j<k\leq n}\left|x_j-x_k\right|^{2/\al}
$$
so we have that the left hand side of (\ref{sgarden}) can be 
rewritten as
\begin{eqnarray*} \lefteqn{
\frac{{\rm C.T.}\left(\bprod_{i=1}^n x_i^{-(n-1)/\al-b}(1-x_i)^{a+b}
E_{\eta}(x_1,\ldots,x_n) \bprod_{j<k}\left|x_j-x_k\right|^{2/\al}
\right) }
{{\rm C.T.}\left(\bprod_{i=1}^n x_i^{-(n-1)/\al-b}(1-x_i)^{a+b}
\bprod_{j<k}\left|x_j-x_k\right|^{2/\al} \right) }  }\hspace{3cm}\\
&&=\frac{\bprod_{l=1}^n\bint^1_0 dt_l t_l^{\lambda_1} (1-t_l)^{\lambda_2}
E_{\eta}(t_1,\ldots,t_n) \bprod_{j<k}\left|t_j-t_k\right|^{2/\al} }
{\bprod_{l=1}^n\bint^1_0 dt_l t_l^{\lambda_1} (1-t_l)^{\lambda_2}
\bprod_{j<k}\left|t_j-t_k\right|^{2/\al} }
\end{eqnarray*}
where the last equality follows from Lemma \ref{plem} with
$\epsilon = -(n-1)/\al -b$, $\lambda_1 = -(n-1)/\al -b+1$, 
$\lambda_2 = a+b$. i.e. $b = -(n-1)/\al -\lambda_1+1$, 
$a=\lambda_1+\lambda_2 +(n-1)/\al +1$. Equating the last equality
to (\ref{sgarden}) gives
$$
\frac{\bprod_{l=1}^n\bint^1_0 dt_l t_l^{\lambda_1} (1-t_l)^{\lambda_2}
E_{\eta}(t_1,\ldots,t_n) \prod_{j<k}\left|t_j-t_k\right|^{2/\al} }
{\bprod_{l=1}^n\bint^1_0 dt_l t_l^{\lambda_1} (1-t_l)^{\lambda_2}
\prod_{j<k}\left|t_j-t_k\right|^{2/\al}  }% \hspace{4cm}\\
= E_{\eta}(1^n) \frac{[\lambda_1 +(n-1)/\al +1]^{(\al)}_{\eta^+} }
{[\lambda_1 + \lambda_2 +2(n-1)/\al +2]^{(\al)}_{\eta^+} }
$$
This formula is the generalization of the Macdonald-Kadell-Kaneko 
integration formula, in the form given by Kaneko \cite{kaneko93a}. The
formula of \cite{kaneko93a} can be reclaimed by multiplying both sides
by $d_{\eta}/e_{\eta}$ and summing over the distinct permutations
of $\ka=\eta^+$ using (\ref{abolladura}).

\setcounter{equation}{0}
\section{The Hermite case}

In this section we shall construct an operator $\ph$ 
analogous to the operator $\Phi$ given in (\ref{phi-op}), which also has
a very simple action on non-symmetric Jack polynomials.  The properties 
of $\ph$ underpin many of the results in this section.

Another key ingredient is the type $A$ Dunkl operator given by
\begin{equation}\label{dunk.op}
T_i = \dif{x_i} +\frac{1}{\al}\sum_{p\neq i} \frac{1}{x_i-x_p}
(1-s_{ip}).
\end{equation}
which satisfies the following relations
\begin{eqnarray}
{[}T_i,x_i]  &=& 1 + {1 \over \alpha} \sum_{p \ne i} s_{ip} \hspace{15mm}
T_i\,s_{ip} = s_{ip} T_p \nonumber\\
{[}T_i,x_j]  &=&  - {1 \over \alpha} s_{ij}, \quad i \ne j \qquad
{}[T_i,s_{jp}] = 0 \quad i\neq j,p \label{bolsilla}
\end{eqnarray}
Note that we can write the Cherednik operator (\ref{cherednik.1})
in the simple form
\begin{equation}\label{simf}
\xi_i = \al x_iT_i + 1-n + \sum_{p>i} s_{ip}
\end{equation}
{}From this and (\ref{bolsilla}) we have the following lemma 
\cite[Lemma 3.1]{forr96c} which will be required later on.
\begin{lemma}\label{desist}
\begin{eqnarray*}
[\xi_j, T_i] & = & T_i s_{ij}, \qquad i < j \\
{}[\xi_j, T_i] & = & T_j s_{ij}, \qquad i > j \\
{}[\xi_j, T_j] & = & -\alpha T_j - \sum_{p < j} s_{jp} T_j - \sum_{p> j} T_j
s_{jp}.
\end{eqnarray*}
\end{lemma}

In \cite{forr96c} we showed that the operators 
\begin{equation}
h_i = \xi_i -\frac{\al}{2}\,T_i^2
\end{equation}
are eigenoperators for the non-symmetric Hermite polynomials $\eh_{\eta}$
and are mutually commuting and self-adjoint with respect to the
inner product (\ref{inn.h}). It is straightforward to show that
the operators $h_i$, $s_j$ generate the same algebra as the operators
$\xi_i$, $s_j$, namely
$$
h_is_i - s_i h_{i+1} = 1, \qquad h_{i+1}s_i -s_i h_i = -1, \qquad
{}[h_i,s_j] = 0, \quad j\neq i, i+1
$$
As a consequence, one can follow the argument of \cite[Prop. 4.3]{knop96c}
and show that 
\begin{equation} \label{sacapuntas}
s_i\,\eh_{\eta} = \left\{ \begin{array}{ll}
\bfrac{1}{\de_{i,\eta}}\,\eh_{\eta} + \left( 1- \bfrac{1}{\de_{i,\eta}^2}
\right) \eh_{s_i\eta} & \eta_i > \eta_{i+1} \\
\eh_{\eta} & \eta_i = \eta_{i+1} \\
\bfrac{1}{\de_{i,\eta}}\,\eh_{\eta} + \eh_{s_i\eta} & \eta_i < \eta_{i+1}
\end{array} \right.
\end{equation}
a result we shall use later on.

\subsection{The operator $\hat{\Phi}$}
We now define the operator $\ph$ as
\begin{equation} \label{escombros}
\ph = T_1 s_1 s_2 \cdots s_{n-1} = s_1s_2\cdots s_{i-1}\; T_i\;
s_i s_{i+1}\cdots s_{n-1}
\end{equation}
This operator obeys the following important relations

\begin{lemma} \label{unfettered}
\begin{eqnarray*}
(a)\qquad \xi_j \;\ph &=& \ph\; \xi_{j-1} \hspace{2cm} \mbox{for}\quad 
2\leq j \leq n \\
(b)\qquad \xi_1\;\ph &=& \ph\;(\xi_n-\al)
\end{eqnarray*}
\end{lemma}
{\it Proof.}\quad First consider (a). From (\ref{hecke}) we have
for $j\geq 2$
\begin{eqnarray*}
\xi_j\;s_1 s_2\cdots s_{n-1} &=& s_1s_2\cdots s_{j-2}\;\xi_j\;
s_{j-1}s_j\cdots s_{n-1} \\
&=& s_1s_2\cdots s_{j-2}\;(s_{j-1}\xi_{j-1} -1)\;s_j\cdots s_{n-1} \\
&=& s_1s_2\cdots s_{n-1}\xi_{j-1} - s_1s_2\cdots s_{j-2}\;
s_j\cdots s_{n-1}
\end{eqnarray*}
{}From this and Lemma \ref{desist} we thus have
\begin{eqnarray*}
\xi_j\;\ph &=& \xi_j\;T_1\;s_1 s_2\cdots s_{n-1} = T_1(\xi_j+s_{1j})
s_1 s_2\cdots s_{n-1} \\
&=& T_1\left(s_1s_2\cdots s_{n-1}\xi_{j-1} - s_1s_2\cdots s_{j-2}\;
s_j\cdots s_{n-1} + s_{1j}\:s_1s_2\cdots s_{n-1} \right)
\end{eqnarray*}
But the permutations occuring in the last two terms in the above
equation are equal, both being equal to $(1\;2\;\ldots\;j-1)
\:(j\;j+1\;\ldots\;n)$ in cycle notation. The result now follows.

Turning to (b), repeated use of (\ref{hecke}) yields
$$
\xi_1\;s_1 s_2\cdots s_{n-1} = s_1 s_2\cdots s_{n-1}\;\xi_n +
\sum_{j=1}^{n-1} s_1\cdots s_{j-1}\;s_{j+1}\cdots s_{n-1}
$$
Thus use of Lemma \ref{desist} gives
\begin{eqnarray*}
\xi_1\;\ph &=& T_1\left( \xi_1 -\al - \sum_{p>1} s_{1p} \right)
s_1s_2\cdots s_{n-1} \\
&=& T_1\left( s_1 s_2\cdots s_{n-1}\:(\xi_n-\al) +
\sum_{j=1}^{n-1} s_1\cdots s_{j-1}\;s_{j+1}\cdots s_{n-1}
- \sum_{p>1} s_{1p}\;s_1s_2\cdots s_{n-1} \right) \\
&=& \ph(\xi_n-\al)
\end{eqnarray*}
where the last equality follows since again the permutations in the
line above cancel.
\hfill$\Box$

\begin{cor} \label{antemano}
The action of $\ph$ on non-symmetric Jack polynomials is given by
$$
\ph\;E_{\eta} = \frac{1}{\al}\;\frac{d'_{\eta}}{d'_{\ph\eta}} \;
E_{\ph\eta}
$$
where $\ph\eta := (\eta_n -1,\eta_1,\eta_2,\ldots,\eta_{n-1})$.
\end{cor}
{\it Proof.}\quad The previous lemma implies that $\ph\;E_{\eta}$
is a constant multiple of $E_{\ph\eta}$. To determine this constant,
note that the leading term in $E_{\ph\eta}$, and hence in
$\ph\;E_{\eta}$ is a multiple of $x^{\ph\eta}$. Writing
$\ph=s_1\cdots s_{n-1} T_n$, we see that the coefficient of $x^{\ph\eta}$
in $\ph\,E_{\eta}$ is equal to the coefficient of 
$x_1^{\eta_1}\cdots x_{n-1}^{\eta_{n-1}}x_n^{\eta_n-1}$ in 
$T_n\,E_{\eta}$. Recalling that $\xi_n = \al x_n\,T_n +1-n$, we
finally deduce that the coefficient of $x^{\ph\eta}$
in $\ph\,E_{\eta}$ is just $(\bar{\eta}_n +n-1)/\al$. Thus
$$
\ph\;E_{\eta} = \frac{\bar{\eta}_n +n-1}{\al}\;E_{\ph\eta}
$$
However, a simple change of variables in Lemma \ref{recur.1}
(recall (\ref{necio})) tells us that
\begin{equation}\label{jude}
\frac{d'_{\eta}}{d'_{\ph\eta}} = \bar{\eta}_n +n-1
\end{equation}
whence the result. \hfill $\Box$

The raising (resp. lowering) operator $\Phi$ (resp. $\ph$) for the
non-symmetric Jack polynomials have their counterparts for the
non-symmetric Hermite polynomials. In fact $\ph$ remains a lowering
operator for the $\eh_{\eta}$, but $\Phi$ no longer has a simple action
in the Hermite case. We find instead that $\ph^{\ast}$ is the
appropriate raising operator for the $\eh_{\eta}$'s, where ${}^{\ast}$
denotes the adjoint operator with respect to the Hermite inner product
(\ref{inn.h}). To show how this comes about, we need some preliminary
results from \cite{forr96c,dunkl89a,dunkl91a}. Following Dunkl, define 
$$
\Delta_A := \sum_{i=1}^n T_i^2
$$
Then we have the commutation relations
\begin{equation}\label{dedo.1}
{}[\xi_i,\Delta_A] = -2\al T_i^2, \hspace{3cm}
{}[x_i,\Delta_A] = -2 T_i .
\end{equation}
Also, the adjoint of the Dunkl operator $T_i$ under the Hermite inner
product (\ref{inn.h}) is given by
\begin{equation}\label{dedo.2}
T_i^{\ast} = 2x_i - T_i
\end{equation}
Finally, using (\ref{dedo.1})  we have the fact that
\begin{equation}\label{dedo.3}
{}[\Phi,\Delta_A] = [s_{n-1}\cdots s_1\;x_1, \Delta_A]
=s_{n-1}\cdots s_1\;[x_1,\Delta_A] = -2s_{n-1}\cdots s_1\;T_1
\end{equation}
The identities (\ref{dedo.2}) and (\ref{dedo.3}) allow a convenient
representation of the operator $\ph^{\ast}$, namely
\begin{eqnarray}
\ph^{\ast} &=& s_{n-1}\cdots s_1\;T_1^{\ast} = 
s_{n-1}\cdots s_1\;(2x_1-T_1) \nonumber\\
&=& 2\Phi + \sfrac{1}{2}[\Phi,\Delta_A] \label{rem}
\end{eqnarray}
We are now in a position to state and prove the Hermite analogues of 
Lemmas \ref{sahi.2} and \ref{unfettered} which take the form

\begin{prop} \label{acotar}
The operators $h_i$ satisfy
\begin{eqnarray*}
(a)\qquad h_n\;\ph^{\ast} &=& \ph^{\ast}\;(h_1+\al) \hspace{2cm}
h_i\;\ph^{\ast} = \ph^{\ast}\;h_{i+1} \qquad 1\leq i\leq n-1 \qquad\\
(b)\qquad h_1\;\ph\;&=& \ph\;(h_n-\al) \hspace{2cm}
h_i\;\ph = \ph\;h_{i-1} \qquad 2\leq i\leq n 
\end{eqnarray*}
\end{prop}
{\it Proof.}\quad We prove only (a) as the proof of (b) is 
straightforward. First note that for $1\leq i\leq n-1$ we have 
$T^2_i\;\Phi = \Phi\;T^2_{i+1} - \gamma$ where
$$
\gamma:= -\frac{1}{\al}s_{n-1}\cdots s_1\;
s_{i+1,1}\left( T_{i+1} + T_1\right) .
$$
This fact, along with Lemma \ref{sahi.2}, (\ref{dedo.1}) and
(\ref{rem}) facilitates a simple calculation which gives 
$ \xi_i\;\ph^{\ast} = \ph^{\ast}\;\xi_{i+1} -\al\gamma$. 

Also, the fact that
the Dunkl operators commute amongst themselves, and hence with $\Delta_A$
tells us that $T^2_i\;\ph^{\ast} = \ph^{\ast}\;T^2_{i+1}-2\gamma$. 
Combining these two results yields the second equality in (a).

In a similar manner, one can show that
$$
\xi_n\;\ph^{\ast} = \ph^{\ast}\;(\xi_1+\al) +\al\;s_{n-1}\cdots
s_1\left( 2T_1 + \frac{1}{\al}\sum_{p>1}(s_{1p}T_i + s_{1p}T_p)
\right)
$$
and
$$
T_n^2\;\ph^{\ast} = \ph^{\ast}\;T_1^2 + 2\;s_{n-1}\cdots
s_1\left( 2T_1 + \frac{1}{\al}\sum_{p>1}(s_{1p}T_i + s_{1p}T_p)
\right)
$$
from which the first equality in (a) follows. \hfill $\Box$

\begin{cor} \label{hola}
The operators $\ph$ and $\ph^{\ast}$ act on the non-symmetric Hermite
polynomials as
\begin{eqnarray}
\ph\;\eh_{\eta} &=& \frac{1}{\al}\;\frac{d'_{\eta}}{d'_{\ph\eta}} \;
\eh_{\ph\eta} \label{hola.1}\\
\ph^{\ast}\;\eh_{\eta} &=& 2\;\eh_{\Phi\eta} \label{hola.2}
\end{eqnarray}
\end{cor}
{\it Proof.}\quad Follows from Proposition \ref{acotar} and
examination of the leading terms on
both sides of the equations. \hfill $\Box$

\vspace{2mm}\noindent
We are now in a position to compute the norm $\innh{\eh_{\eta}}
{\eh_{\eta}}$ in the spirit of the calculation done earlier in the
Jack case.

\begin{prop}\label{tierra}
We have
\begin{equation}\label{lluvia}
\innh{\eh_{\eta}}{\eh_{\eta}} = \frac{1}{(2\al)^{|\eta|}}\;
\frac{d_{\eta}'e_{\eta}}{d_{\eta}}\;{\cal N}_0^{(H)}
\end{equation}
where 
$$
{\cal N}_0^{(H)} := \innh{1}{1} = 2^{-n(n-1)/2\alpha}
\pi^{n/2} \prod_{j=0}^{n-1}{\Gamma(1+(j+1)/\alpha) \over
\Gamma (1 + 1/\alpha )}
$$
is the ground state normalization. 
\end{prop}
{\it Proof.}\quad It is clear that the operators $s_i$ and $\ph^{\ast}$
generate all $\eh_{\eta}$ via (\ref{sacapuntas}) and Corollory 
\ref{hola}. Indeed, as in the Jack case, (\ref{sacapuntas}) and the 
orthogonality of the non-symmetric Hermite polynomials show that
in the case $\eta_i < \eta_{i+1}$
\begin{equation}
\innh{E_{s_i\eta}}{E_{s_i\eta}} =  
(1-\de_{i,\eta}^{-2})\:\innh{E_\eta}{E_\eta} \label{hito.3}
\end{equation}
Also, a simple change of variables in (\ref{hola.1}) gives
$$
\ph\;\eh_{\Phi\eta} = \frac{1}{\al}\;\frac{d'_{\Phi\eta}}{d'_{\eta}} \;
\eh_{\eta}
$$
Hence taking the inner product of both sides of (\ref{hola.2}) with
$\eh_{\Phi\eta}$ and dividing by 2 gives
\begin{eqnarray}
\innh{\eh_{\Phi\eta}}{\eh_{\Phi\eta}} &=& \frac{1}{2}\innh{\eh_{\Phi\eta}}
{\ph^{\ast}\eh_{\eta}} = \frac{1}{2}\innh{\ph\eh_{\Phi\eta}}{\eh_{\eta}}
\nonumber\\
&=& \frac{1}{2\al}\;\frac{d'_{\Phi\eta}}{d'_{\eta}}
\innh{\eh_{\eta}}{\eh_{\eta}}\label{hito.4}
\end{eqnarray}
It thus suffices to show that the right hand side of (\ref{lluvia})
satisfy the recursions (\ref{hito.3}) and (\ref{hito.4}), and is valid
in the case $\eta=0$. But the recursions 
follow immediately from Lemma \ref{recur.1}, while the $\eta=0$ case
is a well-known limiting case of Selberg's integral. \hfill $\Box$

\subsection{Integral kernel and generating function}

The operator $\ph$ has another application concerning integral
kernels introduced by Dunkl \cite{dunkl91a,dunkl92a}. In particular
we are able to derive a generating function for the non-symmetric
Hermite polynomials, which can be used to derive an integral
transform which makes explicit the theory of Dunkl in the type $A$
case. We begin with the following lemma

\begin{lemma} \label{symm}
Let $F=\sum_{\eta} A_{\eta}\,E_{\eta}(x) E_{\eta}(y)$. Then 
$s^{(x)}_i\,F = s^{(y)}_i\,F$ if and only if the coefficients $A_{\eta}$
satisfy
$$
A_{s_i\eta} = \left\{ \begin{array}{cc}
\left( 1- \bfrac{1}{\de_{i,\eta}^2} \right) A_{\eta} & \eta_i > \eta_{i+1} \\
\left( 1- \bfrac{1}{\de_{i,\eta}^2} \right)^{-1} A_{\eta} & \eta_i 
< \eta_{i+1} \end{array} \right.
$$
Moreover, these two conditions on $A_{\eta}$ are equivalent.
\end{lemma}
{\it Proof.}\quad
For a given $i$, split the sum occuring in $s_i^{(x)}\,F$ according to
whether $\eta_i > \eta_{i+1}$, $\eta_i = \eta_{i+1}$ or
$\eta_i < \eta_{i+1}$. Apply Lemma \ref{sahi.1} and collect coefficients 
of $E_{\eta}(x)$. The resulting terms can be identified with $s_i^{(y)}\,F$
if and only if the conditions of the Lemma are satisfied. The
equivalence of the two stated conditions follows through the change
of variable $\eta\rightarrow s_i\eta$. \hfill $\Box$

Let us now define the function
\begin{equation} \label{ker-A}
{\cal K}_A(x;y) = \sum_{\eta} \al^{|\eta|}\;\frac{d_{\eta}}
{d'_{\eta}e_{\eta}}\;E_{\eta}(x)\:E_{\eta}(y)
\end{equation}
Its fundamental properties are given by the following result

\begin{thm} \label{main}
The function ${\cal K}_A(x;y)$ possesses the following properties
\begin{eqnarray*}
(a)\qquad s_i^{(y)}\;{\cal K}_A(x;y) &=& s_i^{(x)}\;{\cal K}_A(x;y) \\
(b)\qquad \ph^{(y)}\;{\cal K}_A(x;y) &=& \Phi^{(x)}\;{\cal K}_A(x;y) \\
(c)\qquad T_i^{(y)}\;{\cal K}_A(x;y) &=& x_i\;{\cal K}_A(x;y) 
\end{eqnarray*}
\end{thm}
{\it Proof.}\\
(a)\quad From Lemma \ref{recur.1}, the constants
$A_{\eta} = \al^{|\eta|}\;\bfrac{d_{\eta}}{d'_{\eta}e_{\eta}}$ 
satisfy the conditions of Lemma \ref{symm} and so the result follows.\\
(b) Using Lemma \ref{sahi.2}, Lemma \ref{recur.1} and Corollory 
\ref{antemano} we have
\begin{eqnarray*}
\ph^{(y)}\;{\cal K}_A(x;y) &=& \sum_{\eta} \al^{|\eta|}\;\frac{d_{\eta}}
{d'_{\eta}e_{\eta}}\;E_{\eta}(x)\; \frac{1}{\al}\;
\frac{d'_{\eta}}{d'_{\ph\eta}} \; E_{\ph\eta}(y) \\
&=& \sum_{\nu} \al^{|\nu|} \frac{d_{\Phi\nu}}{d'_{\nu}e_{\Phi\nu}}\;
E_{\Phi\nu}(x)\:E_{\nu}(y) \\
&=& \sum_{\nu} \al^{|\nu|} \frac{d_{\nu}}{d'_{\nu}e_{\nu}}\;
\Phi^{(x)}\;E_{\nu}(x)\:E_{\nu}(y) \\
&=& \Phi^{(x)}\;{\cal K}_A(x;y)
\end{eqnarray*}
(c)\quad From (\ref{phi-op}) we have
$$
x_i = s_i^{(x)} s_{i+1}^{(x)}\cdots s_{n-1}^{(x)} \;\Phi^{(x)}\;
s_1^{(x)} s_2^{(x)} \cdots s_{i-1}^{(x)}
$$
while from (\ref{escombros}) we have
$$
T_i^{(y)} = s_{i-1}^{(y)} s_{i-2}^{(y)}\cdots s_1^{(y)} \;\ph^{(y)}\;
s_{n-1}^{(y)} s_{n-2}^{(y)} \cdots s_i^{(y)}
$$
Hence using (a), (b) and the fact that operators acting on different
sets of variables commute, we have
\begin{eqnarray*}
T_i^{(y)} {\cal K}_A(x;y) &=& s_{i-1}^{(y)} s_{i-2}^{(y)}\cdots s_1^{(y)}
\;\ph^{(y)}\; s_{n-1}^{(y)} s_{n-2}^{(y)} \cdots s_i^{(y)}\;
{\cal K}_A(x;y) \\
&=& s_i^{(x)}s_{i+1}^{(x)}\cdots s_{n-1}^{(x)}\;
s_{i-1}^{(y)} s_{i-2}^{(y)}\cdots s_1^{(y)} \;\ph^{(y)}\;{\cal K}_A(x;y)\\
&=& s_i^{(x)}s_{i+1}^{(x)}\cdots s_{n-1}^{(x)}\;\Phi^{(x)}\;
s_{i-1}^{(y)} s_{i-2}^{(y)}\cdots s_1^{(y)} \;{\cal K}_A(x;y)\\
&=& s_i^{(x)}s_{i+1}^{(x)}\cdots s_{n-1}^{(x)}\;\Phi^{(x)}\;
s_{1}^{(x)} s_{2}^{(x)}\cdots s_{i-1}^{(x)} \;{\cal K}_A(x;y)\\
&=& x_i\; {\cal K}_A(x;y) 
\end{eqnarray*}
\hfill $\Box$

% The conditions (a)--(c) in Theorem \ref{main} are precisely those
% stipulated by Dunkl in \cite{dunk
The above Theorem enables a generating function for $\eh_{\eta}$ to
be derived. First introduce the notation
\begin{equation}\label{ep}
\et_k:=\sum_{i=1}^n x_i^k \dif{x_i} \hspace{2cm}
p_k(x) := \sum_{i=1}^n x_i^k \;.
\end{equation}
Note that the Hamiltonian (\ref{h-herm}) can be written as
${H}^{(H)} = \Delta_A -2\et_1$. From Theorem \ref{main} (c) we have
\begin{equation}\label{m1}
\Delta_A^{(x)} \;{\cal K}_A(2x;z) = 4\;p_2(z)\;{\cal K}_A(2x;z) .
\end{equation}
We also have that 
\begin{equation}
\et_1^{(z)} \, E_\eta (z) = |\eta|\,  E_\eta (z) \quad
\mbox{and} \quad 
{\cal K}_A(2x;z)\,e^{-p_2(z)} = 
\sum_{\eta} {(2 \alpha)^{|\eta|} d_{\eta} \over
d_\eta' e_\eta} Q_\eta (x) E_\eta (z),
\label{s1}
\end{equation}
where $Q_\eta(x)$ is a polynomial with leading term $E_\eta(x)$.
Following the proof of \cite[Prop.~3.1]{forr96a}, which is Lassalle's
\cite{lass96a}
derivation of the generating function for
symmetric Hermite polynomials, we see by applying the
operator $\et_1^{(z)}$ to both sides of (\ref{ker-A}) and using 
(\ref{s1}) that $Q_\eta(x)$ is an eigenfunction of $H^{(H)}$ with
eigenvalue $-2 |\eta|$. Since the leading term of  $Q_\eta(x)$ is
 $E_\eta(x)$ it follows that $Q_\eta(x) = E^{(H)}_\eta(x)$, and
 thus the generating function for $ E^{(H)}_\eta(x)$ is
 given by
\begin{prop}\label{hermgf}
We have
$$
\sum_{\eta} \frac{(2\al)^{|\eta|}\,d_{\eta}}{e_{\eta}d'_{\eta}}\;
\eh_{\eta}(x)\; E_{\eta}(z) = {\cal K}_A(2x;z)\; e^{-p_2(z)}.
$$
\end{prop}

\vspace{2mm}
An immediate application of Proposition \ref{hermgf} is to provide
an alternative derivation of the norm (\ref{lluvia}). This also
requires using the orthogonality of $\{ E^{(H)}_\eta \}$ with
respect to (\ref{inn.h}) (see the proof of 
\cite[Prop.~3.7]{forr96a} for details, where the
analogous calculation is presented in the symmetric case).

We note that the analytic nature of ${\cal K}_A(x,y)$ is easily 
established.
\begin{prop}
The function ${\cal K}_A(x;y)$ is an entire function of all variables. 
\end{prop}

\vspace{.2cm}
\noindent
{\it Proof.}\quad Since ${\cal K}_A(x;y)$ is a sum of analytic
functions (polynomials), it is sufficient to show that 
$|{\cal K}_A(x;y)|$ is bounded. Now, since the coefficients in the
monomial expansion of
$E_\eta(x)$ are positive \cite{knop96c}, 
$$
|E_\eta(x)| \le E_\eta(1^n) X^{\eta} = {e_\eta \over d_\eta} X^{\eta}
$$
where $X = \mbox{max}(|x_1|, \dots, |x_n|)$, and similarly for 
$|E_\eta(y)|$. Thus
$$
|{\cal K}_A(x;y)| \le \sum_\eta |\alpha XY|^{|\eta|} {1 \over d_\eta'}
E_\eta(1^n) = \sum_\kappa |\alpha XY|^{|\kappa|} {1 \over j_\kappa}
J_\kappa(1^n),
$$
where to obtain the last equality we have used the formula
(\ref{abolladura}). But in general \cite{stan89a}
$$
\sum_\kappa \alpha^{|\kappa|}  {1 \over j_\kappa} J_\kappa(x)
= e^{p_1(x)}
$$
so we obtain the bound
\begin{equation}\label{bound}
|{\cal K}_A(x;y)| \le e^{nXY }
\end{equation}
and the result follows. \hfill$\Box$
\vspace{.2cm}

Next we will show that ${\cal K}_A(x;y)$ is closely related to the
generalized hypergeometric function
$$
{}_0^{} {\cal F}_0^{(\alpha)}(x;y) := \sum_\kappa
{C_\kappa^{(\alpha)}(x) C_\kappa^{(\alpha)}(y)
\over |\kappa|! C_\kappa^{(\alpha)}(1^n)}
$$
where $C_\kappa^{(\alpha)}(x) := \alpha^{|\kappa|} 
|\kappa|! j_\kappa^{-1} J_\kappa^{(\alpha)}(x)$ (below we will
typically omit the superscript $(\alpha)$).

\begin{prop}\label{ss1}
Let {\rm Sym} denote the operation of (full) symmetrization of a function
of $n$ variables, so that
$$
{\rm Sym}\:f(x_1,\dots,x_n) := \sum_P f(x_{P(1)},\dots,x_{P(N)}),
$$
where $P$ denotes a permutation. We have
\begin{equation}
{\rm Sym}^{(x)} {\cal K}_A(x;y) = n! {}_0 {\cal F}_0^{(\alpha)}(x;y).
\label{s5}
\end{equation}
\end{prop}

\vspace{.2cm}
\noindent {\it Proof.}\quad From \cite[eq. (2.18)]{forr96c}, we know that
\begin{equation}\label{symm.1}
{\rm Sym}^{(x)} \, E_\eta (x) = a_\eta J_{\eta^+}(x).
\end{equation}
for some constant $a_{\eta}$. Therefore
\begin{eqnarray}
\mbox{Sym}^{(x)} {\cal K}_A(x;y) & = &
\sum_\eta \alpha^{|\eta|} {d_\eta \over d_\eta' e_\eta}
\mbox{Sym}^{(x)} \, E_\eta(x) E_\eta(y) \nonumber \\
& = & \sum_\eta \alpha^{|\eta|} {d_\eta \over d_\eta' e_\eta}
a_\eta J_{{\eta}^+}(x) E_\eta (y),
\label{s7}
\end{eqnarray}
But from eq.~(a) in Theorem \ref{main}
$$
\mbox{Sym}^{(x)} {\cal K}_A(x;y) = \mbox{Sym}^{(y)} {\cal K}_A(x;y)
$$
and so
\begin{eqnarray}&
\mbox{Sym}^{(x)} {\cal K}_A(x;y)  = 
{1 \over n!} \mbox{Sym}^{(y)} \mbox{Sym}^{(x)} \,  {\cal K}_A(x;y) =
{1 \over n!} \sum_{\eta} \alpha^{|\eta|}  {d_\eta \over d_\eta' e_\eta}
a_\eta^2  J_{{\eta}^+}(x)  J_{{\eta}^+}(y)&  \nonumber
\\ & = \sum_\kappa c_\kappa C_\kappa(x) C_\kappa(y)/C_\kappa(1^n) &
\label{s3}
\end{eqnarray}
for some $c_\kappa$.

To determine $c_\kappa$, from eq.~(c) of Theorem \ref{main} and
the fact that $\sum_i T_i^{(y)} = \tilde{E}_0^{(y)}$, we note
$$
 \tilde{E}_0^{(y)} \,  {\cal K}_A(x;y) = p_1(x) \,  {\cal K}_A(x;y)
$$
But $p_1(x)$ is symmetrical in $x$, so after 
applying $\mbox{Sym}^{(x)}$ to both sides and using (\ref{s3}) 
we find 
$$
\sum_{\kappa} c_\kappa  \tilde{E}_0^{(y)} \, C_\kappa(x) C_\kappa(y)
= \sum_{\kappa} c_\kappa p_1(x)  C_\kappa(x) C_\kappa(y).
$$
Using known formulas for $ \tilde{E}_0^{(y)} \,  C_\kappa(y)$ and
$ p_1(x)  C_\kappa(x)$ \cite[eqs.~(2.13a),(3.4)]{forr96a} we see
that this equation uniquely determines the $c_\kappa$ as
$c_\kappa = c_0/|\kappa|!$. But from the definition of Sym we have,
Sym 1 = n!, and so $c_0 = n!$. 
\hfill $\Box$

\vspace{.2cm}
A spin-off from Proposition \ref{ss1} is the explicit value of the
constant $a_\eta$ in (\ref{symm.1}).
\begin{prop}
Let $a_\eta$ be defined by (\ref{symm.1}). We have
\begin{equation}\label{symm.2}
a_\eta = n! {e_\eta \over d_\eta} {1 \over J_\kappa(1^n)}.
\end{equation}
\end{prop}

\vspace{.2cm}
\noindent
{\it Proof.} \quad
We first note from (\ref{abolladura}) that
$$
{C_\kappa(y) \over C_\kappa (1^n)} = {J_\kappa (y) \over J_\kappa(1^n)}
= {j_\kappa \over J_\kappa(1^n)} \sum_{\eta: \, \eta^+ = \kappa}
{1 \over d_{\eta'}} E_\eta(y).
$$
Substituting in the right hand side of
(\ref{s5}) and equating coefficients of
$J_\kappa (x) E_\eta (y)$ with the right hand side of
(\ref{s7}) multiplied by $n!$ gives the sought result.\hfill$\Box$

\vspace{.2cm}
The hypergeometric function ${}_0 {\cal F}_0(x;y)$ is
related to the symmetric Hermite polynomials by a generating function
analogous to Proposition \ref{hermgf}, and also through an integral
transform of the symmetric Jack polynomials, in which 
${}_0 {\cal F}_0(x;y)$ is the kernel \cite{lass96a,forr96a}.
Likewise, ${\cal K}_A(x;y)$ also occurs as the kernel in an
integral transform which relates the non-symmetric Jack
and Hermite polynomials. 

\begin{prop} \label{int.h}
Let 
\begin{equation}
d\mu^{(H)}(y) := \prod_{j=1}^n e^{-y_j^2} \prod_{1 \le j < k \le n}
|y_j - y_k|^{2/\alpha} \, dy_1 \dots dy_n.\label{hm}
\end{equation}
Then we have
\begin{eqnarray}
\int_{(-\infty,\infty)^n} {\cal K}_A(2y;z) {\cal K}_A(2y;w)\: d\mu^{(H)}(y) &=&
{\cal N}_0^{(H)}\,e^{p_2(w)+p_2(z)}\,{\cal K}_A(2z;w) \label{1a} \\
\int_{(-\infty,\infty)^n} {\cal K}_A(2y;z)\;\eh_{\eta}(y)\: d\mu^{(H)}(y) &=&
{\cal N}_0^{(H)}\,e^{p_2(z)}\;E_{\eta}(z) \label{1b} \\
\int_{(-\infty,\infty)^n} {\cal K}_A(2y;-iz)\;E_{\eta}(iy)\: d\mu^{(H)}(y) &=&
{\cal N}_0^{(H)}\,e^{-p_2(z)}\;\eh_{\eta}(z) \label{1c}
\end{eqnarray}
\end{prop}

\vspace{.2cm}
\noindent
{\it Proof.}\quad We first note from the bound (\ref{bound}) that all the
integrals exist. The first formula follows by multiplying both sides 
by $e^{-p_2(w)-p_2(z)}$, using the generating function of Proposition
\ref{hermgf} twice on the left hand side and then using the
orthogonality of $\{ E_\eta^{(H)}(y) \}$ with respect to (\ref{inn.h})
to compute the integral. The resulting sum is identified as
${\cal K}_A(2z;w)$. The second formula follows from the first after
multiplying by $e^{-p_2(w)}$, using the generating function on the
left hand side and equating coefficients of $E_\eta(w)$ on both sides,
while the third follows by replacing $z$ by $iz$, using the
generating function on the right hand side and equating coefficients
of $E_\eta(w)$.
\hfill $\Box$
\vspace{.2cm}

There are a number of consequences of Proposition \ref{int.h}. First
we note a summation theorem, which is the non-symmetric analogue of
\cite[Proposition 3.9]{forr96a}.

\begin{prop} \label{green.h}
For $|t| < 1$ we have
\begin{eqnarray*}
\sum_{\eta} \frac{1}{{\cal N}_{\eta}^{(H)}}\;
\eh_{\eta}(w)\:\eh_{\eta}(z) t^{|\eta|} 
= \frac{1}{{\cal N}_0^{(H)}} (1-t^2)^{-nq/2} \hspace{4cm}\\
\times\exp\left(-\frac{t^2}{(1-t^2)}\left(p_2(z) + p_2(w)\right)\right)\;
{\cal K}_A\left(
\frac{2wt}{(1-t^2)^{1/2}};\frac{z}{(1-t^2)^{1/2}}\right )
\end{eqnarray*}
where $q:=1+(n-1)/\al$.
\end{prop}

\vspace{.2cm}
\noindent
{\it Proof.} \quad This follows by substituting the integral
representation (\ref{1a}) in the left hand side, as in the proof
of \cite[Proposition 3.9]{forr96a}.
\hfill $\Box$
\vspace{.2cm}

The sum in Proposition \ref{green.h} is closely related to the Green
function solution  of the imaginary time Scr\"odinger equation
\begin{equation}
\Big ( \psi_0^{(H)} H^{(H)} (\psi_0^{(H)})^{-1} \Big )G =
{\partial \over \partial \tau} G, \quad
 \psi_0^{(H)} :=  \prod_{l=1}^n 
 e^{-x_l^2} \prod_{1 \le j < k \le n} |x_k - x_j|^{2/\alpha},
 \end{equation}
 which is the solution subject to the initial condition
 $$
 G(x^{(0)}|x;\tau) \Big |_{\tau = 0} = \prod_{j=1}^n \delta(x_j-x_j^{(0)}).
 $$
 Indeed, since $\{E_\eta^{(H)} \}$ is a complete set of eigenfunctions
 of $H^{(H)}$ which are orthogonal with respect to the inner product
 (\ref{inn.h}), by the method of separation of variables we have
\begin{equation} \label{gr.h} 
 G(x^{(0)}|x;\tau) = \psi_0^{(H)}(x^{(0)}) \psi_0^{(H)}(x)
 \sum_{\eta} {1 \over {\cal N}_\eta^{(H)}} \eh_{\eta}(x^{(0)})
 \eh_{\eta}(x) e^{-2\tau |\eta|}.
\end{equation} 
 We can use (\ref{gr.h}) to determine the large $x$ (or large $y$)
 asymptotic behaviour of ${\cal K}_A(x;y)$, since for
 $\tau \to 0$ the asymptotic form of $ G(x^{(0)}|x;\tau)$
must agree with the Green function solution of
 $$
 \sum_j {\partial^2 \over \partial x_j^2} G = 
 {\partial \over \partial \tau} G,
 $$
 which gives
 $$
  G(x^{(0)}|x;\tau) \: \sim \: \Big (
  {1 \over 4 \pi \tau} \Big )^{n/2} \prod_{j=1}^n
  e^{-(x_j - x_j^{(0)})^2/4 \tau}.
$$
Substituting in (\ref{gr.h}) and use of Proposition \ref{green.h}
shows
\begin{prop}
We have
\begin{equation}\label{3.30}
{\cal K}_A(x/\tau^{1/2};y/\tau^{1/2}) \: \sim \:
{\pi^{-n/2} 2^{n(n-1)/2\alpha}  {\cal N}_0^{(H)} \over
( \prod_{1 \le j < k \le n} (x_j - x_k)(y_j - y_k)/\tau)^{1/\alpha}}
 \prod_{j=1}^n e^{x_j y_j/\tau}.
\end{equation}
\end{prop}

We note that (\ref{3.30})
is identical to the asymptotic behaviour of ${}_0{\cal F}_0(x/\tau
;y/\tau)$ (with the assumption $x_1 < \cdots < x_n$ and
$y_1 < \cdots y_n$) \cite[eq.~(5.46)]{forr96a}.

Dunkl \cite{dunkl91a,dunkl92a}
has developed a theory of integral transforms within the framework
of root systems. In the type A case, the kernel satisfies properties
(a) and (b) of Theorem 3.8. We can see that the type A
kernel of Dunkl must be precisely our ${\cal K}_A(x,y)$ by further
developing the consequences of Proposition \ref{int.h}. For this we
require the exponential operator formula \cite{forr96c,sogo96}
\begin{equation}\label{exp.g}
\eh_{\eta} = e^{-\Delta_A/4}\;E_{\eta}
\end{equation}
(an alternative derivation of this formula is afforded by using
(\ref{m1}) to deduce that $e^{-\Delta_A^{(x)}/4}
{\cal K}_A (2x;z) = e^{-p_2(z)} {\cal K}_A (2x;z)$ and then using the
generating function to equate coefficients of $E_\eta^{(H)}(z)$
on both sides). Substituting (\ref{exp.g}) in (\ref{1b}) and
(\ref{1c}), and using the fact that $\{E_\eta^{(H)}\}$ is a basis
for analytic functions, we obtain the following formulas.

\begin{prop}\label{3.zz}
We have
\begin{equation}\label{2a}
\int_{(-\infty,\infty)^n} {\cal K}_A(2y;z)
\Big (  e^{-\Delta_A/4}\;f(y) \Big ) \,  d\mu^{(H)}(y) =
 {\cal N}_0^{(H)} e^{p_2(z)} f(z)
 \end{equation}
\begin{equation}\label{2b}
\int_{(-\infty,\infty)^n} {\cal K}_A(2y;z)\;f(iy) \: d\mu^{(H)}(y)
= {\cal N}_0^{(H)} e^{-p_2(z)}  e^{-\Delta_A/4}\; f(z)
\end{equation}
where $f$ is an analytic function such that all terms converge.
\end{prop}

Now the type A kernel of Dunkl, here denoted $K_A(x,y)$, has the property
(\ref{2a}) with ${\cal K}_A$ replaced by $K_A$
\cite[Prop.~2.1]{dunkl91a}. Since $f$ is arbitrary we must
therefore have
\begin{equation}
K_A(x;y) = {\cal K}_A(x;y).
\end{equation}

Continuing the development of the implications of Proposition \ref{3.zz} 
we note from (\ref{2b}) that if we define an integral
transform (a generalized Fourier transform) by
\begin{equation}\label{gFt}
F(z)  =  {e^{p_2(z)} \over  {\cal N}_0^{(H)}}
\int_{(-\infty,\infty)^n} {\cal K}_A(2y;z)\;f(iy) \: d\mu^{(H)}(y),
\end{equation}
where it is assumed that the integral is absolutely convergent,
then $F$ is inverted by
\begin{equation}\label{inv.d}
f(z) = \exp \Big ( {1 \over 4} \Delta_A \Big ) F(z).
\end{equation}
To obtain the inversion as an integral transform, we follow the 
method given in ref.~\cite{forr96a} for the symmetric case, and
replace $z$ by $iz$ and $f(ix)$ by $F(x)$ in (\ref{2b}),
which gives
\begin{equation}\label{gFt.i}
f(z) = {e^{-p_2(z)} \over  {\cal N}_0^{(H)}}
\int_{(-\infty,\infty)^n}  {\cal K}_A(2y;z)\; F(y) \:  d\mu^{(H)}(y)
\end{equation}
Comparison with (\ref{inv.d}) shows

\begin{prop} \label{inv.ft}
Let $F$ be given in terms of $f$ by (\ref{gFt}). Then
$f$ is given in terms of $F$ by
\begin{equation}
f(z) = \frac{1}{{\cal N}_0^{(H)}}\;e^{-p_2(z)}
\int_{(-\infty,\infty)^n} {\cal K}_A(2y;z)\;F(y) \: d\mu^{(H)}(y)
\end{equation}
\end{prop}

There is further overlap with the theory of Dunkl. For homogeneous
polynomials $p$ and $q$ of the same degree $|\eta|$ say, let
\begin{equation}\label{pairing}
[p,q]_H = p(T^x) q(x),
\end{equation}
where $ p(T^x)$ means each variable $x_i$ in $p$ is replaced by the
Dunkl operator $T_i$. If $p$ and $q$ have different degrees set
$[p,q]_H = 0$. The theory of Dunkl \cite[Thm.~3.10]{dunkl91a} 
gives that
$$
[p,q]_H = {2^{|\eta|} \over {\cal N}_0^{(H)}} \int_{(-\infty,\infty)^n}
\Big ( e^{-\Delta_A/4} p \Big ) \Big ( e^{-\Delta_A/4} q \Big )
\, d \mu^{(H)}(x).
$$
{}From (\ref{exp.g}), (\ref{lluvia}) and the orthogonality of
$\{ E_\eta^{(H)}\}$ with respect to (\ref{inn.h}) we obtain
\begin{prop}
We have
\begin{equation}\label{pairing.h}
[E_\nu,E_\eta]_H = {1 \over \alpha^{|\eta|}} 
{d_\eta' e_\eta \over d_\eta} \delta_{\nu,\eta}.
\end{equation}
\end{prop}

\quad We will conclude this subsection by presenting results, communicated
to us by 
C.~Dunkl \cite{dunkl96c}, on the relationship between the pairing
(\ref{pairing}) and the analogue of the power sum inner product for the 
non-symmetric Jack polynomials. The latter inner product is defined by
Sahi \cite{sahi96a} according to $\langle x^\eta, p_\nu \rangle =
\delta_{\eta,\nu}$, where the polynomials $p_\eta$ are defined by
the expansion 
\begin{equation}\label{dgf}
\prod_{i=1}^n {1 \over 1 - x_i y_i}
\prod_{i,j=1}^n {1 \over (1 - x_i y_j)^{1/\alpha}} =
\sum_{\eta} p_\eta(x) y^\eta
\end{equation}
(Sahi uses the notation $q_\eta$ for $p_\eta$; we use the latter notation
for consistency with ref.~\cite{dunkl96a}). Thus if $f$ and $g$ are
homogeneous polynomials of the same degree, and
\begin{equation}\label{fgp}
f(x) = \sum_\eta f_\eta p_\eta(x), \quad
g(x) = \sum_\eta g_\eta p_\eta(x), \quad
p_\eta(x) = \sum_\nu A_{\eta \nu} x^{\nu},
\end{equation}
then
\begin{equation}\label{sip}
\langle f, g \rangle = \sum_{\eta} \sum_{\nu} f_\eta A_{\eta \nu} g_{\nu}.
\end{equation}

As noted in  \cite{sahi96a}, the generating function (\ref{dgf})
occurs in the recent work of Dunkl \cite{dunkl96a}. Also introduced in
\cite{dunkl96a} is the space $A_\lambda$ (Dunkl uses $E_\lambda$ but this
would cause confusion with the notation for the non-symmetric Jack
polynomials) of homogeneous polynomials of degree $|\lambda|$, where
each $f \in A_\lambda$ has the additional property that
\begin{equation}\label{vf}
V \, \xi \, f = {1 \over [n/\alpha + 1]_\lambda^{(\alpha)}}f.
\end{equation}
Here $\xi$ is defined by $\xi p_\eta = {1 \over \eta_1! \dots \eta_n!}
x^\eta$, and $V$ is defined by the intertwining relation $T_i V = V
{\partial \over \partial x_i}$ with the normalization
$V \, 1 = 1$ and $[c]_{\lambda}^{(\alpha)}$ by (\ref{blitz}). 
Another characterization of $A_\lambda$ is
that it is invariant under the action of $T_i x_i$ (\cite[Prop.~3.2]{dunkl96a}),
and so from (\ref{simf}) and Lemma \ref{sahi.1} it follows that
 $A_\lambda$ is equal to Span$\{E_\eta : \eta^+ = \lambda\}$.

 Let us now suppose $f$ and $g$ in
 (\ref{fgp}) are elements of $A_\lambda$. Then according to (\ref{fgp}),
 (\ref{vf}) and the defining properties of $\xi$ and $V$ we have
 \begin{eqnarray}\label{pfga}
 [f,g]_H & = & \sum_{\nu'} f_{\nu'} \sum_\nu A_{\nu' \nu} T^{\nu}\Big (
 \sum_\eta g_\eta p_\eta(x) \Big ) \nonumber \\
 & = &  [n/\alpha + 1]_\lambda^{(\alpha)} 
 \sum_{\nu'} f_{\nu'} \sum_\nu A_{\nu' \nu} T^{\nu}\Big (
  \sum_\eta  V \xi \, g_\eta p_\eta(x) \Big ) \nonumber \\
  & = & [n/\alpha + 1]_\lambda^{(\alpha)} 
   \sum_{\nu'} f_{\nu'} \sum_\nu A_{\nu' \nu} 
 \sum_\eta {1 \over \eta_1! \dots \eta_n!} g_\eta V \Big ( {\partial \over
 \partial x} \Big )^\nu \, x^\eta \nonumber \\
 & = &
  [n/\alpha + 1]_\lambda^{(\alpha)}  \sum_{\nu'} f_{\nu'} \sum_\nu A_{\nu'
  \nu} g_\nu
  \end{eqnarray}
Comparing (\ref{pfga}) with (\ref{sip}) we see that
\begin{equation}\label{dufgh}
 [f,g]_H =  [n/\alpha + 1]_\lambda^{(\alpha)} \langle f, g \rangle,
 \end{equation}
 which upon using the facts \cite{sahi96a} that
 $e_\eta = \alpha^{|\eta|}[n/\alpha + 1]_\eta$ and
 $\langle E_\eta, E_\nu \rangle = (d_\eta' / d_\eta) \delta_{\nu, \eta}$
 is seen to be 
 in precise agreement with (\ref{pairing.h}).

\subsection{Evaluation of $e^{p_1(x)} {\cal K}_A (x;y)$ and generalized
binomial coefficients}
In the theory of the generalized hypergeometric function
${}_0{\cal F}_0(x;y)$ the identity 
\begin{equation}\label{s.pK}
e^{p_1(x)} {}_0{\cal F}_0(x;y) = {}_0{\cal F}_0(x;y+1),
\end{equation}
is an immediate consequence of the identity \cite{kaneko93a,lass90a}
\begin{equation}\label{conn}
e^{p_1(z)} C_\kappa^{(\alpha)}(z) = \sum_\mu \left (
{ \mu \atop \kappa }\right ) {|\kappa|! \over | \mu |!}
C_\mu^{(\alpha)}(z),
\end{equation}
with the generalized binomial coefficients defined by
\begin{equation}\label{symm.b}
{ C_\kappa^{(\alpha)}(1+z) \over  C_\kappa^{(\alpha)}(1^n)} =
\sum_\sigma \left ( {\kappa \atop \sigma} \right )
{ C_\sigma^{(\alpha)}(z) \over  C_\sigma^{(\alpha)}(1^n)}.
\end{equation}
Analogous results hold for the function ${\cal K}_A (x;y)$, although
it is the analogue of (\ref{s.pK}) which is derived directly.

\begin{prop}
We have
\begin{equation}\label{pK}
e^{p_1(x)}  {\cal K}_A(x;y) =  {\cal K}_A(x;y+1).
\end{equation}
\end{prop}

\vspace{.2cm}
\noindent {\it Proof.} \quad Consider the action of $\xi_i^{(y+1)}$ on
$e^{p_1(x)}  {\cal K}_A(x;y)$. Since $\xi_i^{(y+1)} = \xi_i^{(y)} +
\alpha T_i^{(y)}$ we have
\begin{eqnarray}\label{xii}
\xi_i^{(y+1)} \, e^{p_1(x)}  {\cal K}_A(x;y) & = &
 e^{p_1(x)} \Big (  \xi_i^{(y)} + \alpha T_i^{(y)} \Big )  {\cal K}_A(x;y)
 \nonumber \\
 & = &  e^{p_1(x)} \Big (  \xi_i^{(x)} + \alpha x_i \Big )  {\cal K}_A(x;y)
  \nonumber \\
& = &  \xi_i^{(x)} \Big (  e^{p_1(x)}  {\cal K}_A(x;y) \Big )
\end{eqnarray}
But from the definition of $ {\cal K}_A(x;y)$ we must have
$$
e^{p_1(x)}  {\cal K}_A(x;y) = \sum_{\eta}
\alpha^{|\eta|} { d_\eta \over d_{\eta}' e_\eta} U_\eta(y) E_\eta(x)
$$
where $U_\eta(y)$ is a polynomial with leading term $E_\eta(y)$.
Substituting this in (\ref{xii}) and equating coefficients of
$E_\eta(x)$ we see that
$$
\xi_i^{(y+1)}  U_\eta(y) = \bar{\eta}_i  U_\eta(y)
$$
and thus $ U_\eta(y)$ is proportional to $E_\eta(y+1)$. The proportionality
constant is unity since the leading term of $ U_\eta(y)$ is
$E_\eta(y)$.
\hfill $\Box$

\vspace{.2cm}

We can reclaim (\ref{s.pK}) from (\ref{pK}) by symmetrizing both sides
with respect to $x$ using (\ref{s5}). Since (\ref{s.pK}) and (\ref{conn})
are equivalent, this also establishes the latter identity. A
significant feature of (\ref{conn}) is that combined with the fact
that the coefficients in the monomial expansion of $C_\kappa$ are
independent of the number of variables $n$ it immediately implies
the binomial coefficients $ \left ( {\kappa \atop \sigma} \right )$
are independent of $n$. The only other published proof of the
independence of the binomial coefficients on $n$ is given in
ref.~\cite{kaneko93a}, via a long case-by-case check, and the
independence is then used in the derivation of (\ref{conn}).

In the non-symmetric case we can define generalized
binomial coefficients $\left ( {\eta \atop \nu} \right )$, 
$|\eta| \geq |\nu|$, analogous to (\ref{symm.b}) by
\begin{equation}\label{binomial}
{E_\eta(1+z) \over E_\eta(1^n)}  = \sum_\nu \left ( {\eta \atop \nu} \right )
{E_\nu(z) \over E_\nu(1^n)}
\end{equation}
Some immediate properties of these coefficients are
\begin{lemma}\hfill\\
(i) \hspace{3ex} For $|\nu| = |\eta|$, $ \left ( {\eta \atop \nu} \right ) =1$
for $\nu = \eta$ and $ \left ( {\eta \atop \nu} \right ) =0$ otherwise. \\
(ii) \hspace{2ex} Let $(0,\dots,0) = 0$. Then $ \left ( {\eta \atop 0}
\right ) =1$.
\end{lemma}
Moreover we can provide a simple proof that
\begin{prop}
The coefficients $ \left ( {\eta \atop \nu} \right )$, like their 
symmetric counterparts, are independent of the number of variables $n$.
\end{prop}

\vspace{.2cm}
\noindent
{\it Proof.} \quad This is done in an analogous way to the method of
proof of the independence of the symmetric binomial coefficients noted
above. Thus we substitute the definition (\ref{binomial}) in the 
right hand side of 
(\ref{pK}), and equate
coefficients of $E_\eta(y)$ on both sides to obtain the analogue of
(\ref{conn}):
\begin{equation}\label{n.conn}
e^{p_1(x)} E_\eta(x) \alpha^{|\eta|}{1  \over d_\eta'} =
\sum_{\nu} \alpha^{|\nu|} {1 \over d_{\nu}'} 
\left ( {\nu \atop \eta} \right )E_\nu(x),
\end{equation}
The independence of the $ \left ( {\eta \atop \nu} \right )$ on $n$
follows immediately from this identity, and the facts that the
coefficients $a_{\eta \nu}$ in (\ref{1.3}) and the $d_\eta'$ are
independent of $n$. \hfill$\Box$

\vspace{.2cm}
The symmetric generalized binomial coefficients can be expressed in terms of
the non-symmetric ones. Indeed note that by symmetrizing 
(\ref{binomial}) using (\ref{symm.1})
and (\ref{symm.2}), and comparing with 
(\ref{symm.b}), we obtain
\begin{equation}
\sum_{\nu : \nu^+ = \mu}  \left ( {\eta \atop \nu} \right )
=  \left ( {\kappa \atop \mu} \right ),
\end{equation}
while symmetrizing (\ref{n.conn}) and comparing with (\ref{conn})
gives
\begin{equation}\label{bin.a}
{d_\eta' j_\mu  J_\kappa(1^n) \over  j_\kappa J_\mu(1^n)}
\sum_{\nu:\nu^+=\mu} {1 \over d_\nu'} 
  \left ( {\nu \atop \eta } \right ) = 
  \left ( {\mu \atop \kappa} \right).
\end{equation}

For further application of the non-symmetric  binomial coefficients,
let  $\tilde{E}_0$ and  $\tilde{E}_2$ be defined by (\ref{pK}),
and set
\begin{equation}\label{dt1}
\tilde{D}_1 := D_1 - {1 \over \alpha} \sum_{j \ne k} 
{x_j \over (x_j - x_k)^2} (1 - M_{jk})
\end{equation}
where
$$
D_p := \sum_{j=1}^n x_j^p {\partial^2 \over \partial x_j^2}
+ {2 \over \alpha} \sum_{j \ne k} {x_j^p \over x_j - x_k}
{\partial \over \partial x_j}.
$$
The non-symmetric binomial coefficients can be used to compute the
action of $\tilde{E}_0$, $\tilde{E}_2$ and $\tilde{D}_1$ on the $E_\eta$.
\begin{prop}\label{action.b}
We have
\begin{eqnarray}
\tilde{E_0} \, {E_\eta(x) \over E_\eta(1^n)}& =&
 \sum_{\nu: |\nu| = |\eta| - 1} \left ( {\eta \atop \nu} \right )
 {E_\nu(x) \over E_\nu(1^n)} \\
\tilde{E_2} \, E_\eta(x)&  =& {\alpha \over 2} d_\eta'
\sum_{\nu: |\nu| = |\eta| + 1} {1 \over d_\nu'}
 \left ( {\nu \atop \eta} \right ) (\epsilon_\nu - \epsilon_\eta-{2 \over
 \alpha}(N-1))
  E_\nu(x) \\
  \tilde{D}_1  \, {E_\eta(x) \over E_\eta(1^n)}
&  = & {1 \over 2}\sum_{\nu : |\nu| =
  | \eta| - 1}  \left ( {\eta \atop \nu} \right )
  (\epsilon_\eta - \epsilon_\nu) {E_\nu(x) \over E_\nu(1^n)}
 \end{eqnarray}
 where
 \begin{equation}\label{He}
 \epsilon_\eta = \sum_{j=1}^n \Big ( \eta_j^+ (\eta_j^+ - 1)
 + {2 \over \alpha} (N-j) \eta_j^+ \Big ).
 \end{equation}
 \end{prop}

 \vspace{.2cm}
 \noindent
 {\it Proof.} \quad 
The action of $\tilde{E_0}$ follows immediately from the definition
(\ref{binomial}) used to expand $ E_\eta(\epsilon + x)$ in
the formula 
\begin{eqnarray}
\tilde{E_0} \, E_\eta(x) & = & \lim_{\epsilon \to 0}
{1 \over \epsilon} \Big ( E_\eta(\epsilon + x) - E_\eta(x) \Big ) \nonumber
\\ & = & \sum_{\nu: |\nu| = |\eta| - 1} \left ( {\eta \atop \nu} \right )
E_\nu(x).
\end{eqnarray}
To compute the action of $\tilde{E}_2$ we follow the strategy given
in ref.~\cite{muir82} in the symmetric case and apply the operator 
\begin{eqnarray}\label{dt2}
\tilde{D}_2 & := &  D_2 - {1 \over \alpha} \sum_{j \ne k}
{x_j x_k \over (x_j - x_k)^2} (1 - s_{jk}) \nonumber \\
& = & H^{(C)} + {2 \over \alpha} (N-1) \tilde{E}_1,
\end{eqnarray}
which is an eigenoperator for each $E_\eta$ with corresponding eigenvalue
$\epsilon_\eta$ (\ref{He}), to the identity  
\begin{equation}\label{pex}
p_1(x) E_\eta(x) = \alpha d_\eta' \sum_{\nu: |\nu| = |\eta| +1}
{1 \over d_\nu'}  \left ( {\nu \atop \eta} \right ) E_\nu(x),
\end{equation}
which follows immediately from (\ref{n.conn}).
The action of $\tilde{D}_1$ follows from the action of
$\tilde{E}_0$ and $\tilde{D}_2$ and the readily
verified identity $\tilde{D}_1 = {1 \over 2}[\tilde{E}_0, \tilde{D}_2 ]$.
\hfill $\Box$

\vspace{.2cm}
The formulas of Proposition \ref{action.b} will be used in 
Section 4 to derive a partial differential equation satisfied by
a generalization of ${\cal K}_A(x;y)$, which has application in the
derivation of generating functions for the non-symmetric
Laguerre polynomials.

\subsection{The Laplace transform}
The non-symmetric generalized Laplace transform is defined by
\begin{equation}\label{laplace}
{\cal L}[f](t) = \int_{[0,\infty)}{\cal K}_A(-t;x) f(x)
\prod_{1 \le j < k \le n} |x_k - x_j|^{2/\alpha} dx_1 \dots dx_n
\end{equation}
where it is assumed the integral is absolutely convergent
(note from (\ref{3.30}) that for large-$x$, and $x$ and $t$ suitably
ordered
$$
{\cal K}_A(-t;x) = O (e^{-t_1 x_1 - \dots - t_n x_n}).
$$
In the
symmetric case (${\cal K}_A(-t;x)$ replaced by ${}_0{\cal F}_0(-t;x)$)
with $n=2$, this has been considered in some detail by Yan \cite{yan92a}.
We find that the results of Yan all have non-symmetric counterparts.

A very simple example is the shift property \cite[Prop. 4.10]{yan92a}
\begin{equation}
{\cal L}[e^{-p_1(x)} f](t) = {\cal L}[f] (t+1),
\end{equation}
which follows from the definition and (\ref{pK}). A more fundamental
result is that ${\cal L}$ is injective.

\begin{prop} If $\langle f,f\rangle_L < \infty$ and ${\cal L}[f](t) = 0$ for all $t$ then $f=0$ a.e..
\end{prop}

\noindent
{\it Proof} \quad This is shown by adapting the strategy of Yan.
We will make use of the formula
$$
P_\nu(T^t){E_\eta(t) \over E_\eta(1^n)} \Big |_{t=0} = \delta_{\nu,\eta},
\quad \mbox{where} \quad 
P_\nu(x) = {(2 \alpha)^{|\nu|} \over d_\nu'} E_\nu(x),
$$
 which follows from (\ref{pairing}). Applying this formula to
 (\ref{binomial}) gives that for all $\eta$ with $|\nu| \le 
 |\eta|$
$$
P_\nu(T^t) {E_\eta(t) \over E_\eta(1^n)} \Big |_{t=1} = \Big (
{ \eta \atop \nu} \Big ),
$$
and so
$$
P_\nu(T^t) {\cal K}_A(-t;x) \Big |_{t=1} = (-\alpha)^{|\nu|}
{1 \over d_\nu'} E_\nu(x) e^{-p_1(x)}
$$
where the formula (\ref{n.conn}) has also been used. Thus
$$
P_\nu(T^t) {\cal L}[f](t)  \Big |_{t=1} = {(-\alpha)^{|\nu|} \over
 d_\nu'}  \int_{[0,\infty)} e^{-p_1(x)} E_\nu(x) f(x)
 \prod_{1 \le j < k \le n} |x_k -x_j|^{2/\alpha} dx_1 \dots dx_n.
$$
Since we are assuming ${\cal L}[f](t) = 0$, this last expression
must vanish, and  this holds true for all $\nu$. But
$\{E_\nu\}$ is a basis for analytic functions, so it follows
that $f=0$ (this can been seen by forming an appropriate linear
combination of the $E_\nu$ so as to reconstruct $f(t)$ a.e., which
gives $\langle f,f \rangle_L \Big |_{a=0} = 0$.)
 \hfill$\Box$

\vspace{.2cm}
Further properties of the generalized Laplace transform will
be discussed in the next section.

\subsection{Relationship to Dunkl's theory of harmonic polynomials}
Van Diejen \cite{vandiej96b} has shown how the symmetric generalized
Hermite (and Laguerre) polynomials can be written in terms of Dunkl's
generalized spherical polynomials \cite{dunkl88a}. In ref.~\cite{dunkl88a}
Dunkl has extended the theory of \cite{dunkl84a} to the non-symmetric case,
and this allows the considerations of \cite{vandiej96b} to be similarly
extended.

The $A$ type generalized harmonic polynomials, ${\cal Y}^A_{k,l}$ say, are
defined by the equation
\begin{equation}\label{DYA}
\Delta_A^2 \,{\cal Y}^A_{k,l} = 0
\end{equation}
and the conditions that they are homogeneous of degree $k$ and linearly
independent. The label $l$ distinguishes linearly independent solutions
for each value of $k$; with $P_k$ denoting the
space of homogeneous polynomials of degree $k$,
Dunkl \cite{dunkl89a} has shown that there are dim$\, P_k$ $-$ dim$\,
P_{k-2}$
linearly independent solutions. However, in Dunkl's theory the basis for
the label $l$ is left unspecified, and in the equations below $l$ will be
replaced by a dot.

Now, from \cite[Th.~1.11]{dunkl88a} any homogeneous polynomial can be
expressed
in terms of certain harmonic polynomials, which themselves
are specified by the homogeneous polynomial. Applying this formula 
to the non-symmetric Jack polynomials gives
\begin{equation}\label{ey}
E_\eta(x) = \sum_{m=0}^{[|\eta|/2]} r^{2m} {\cal Y}^A_{|\eta| - 2m,\cdot}(x),
\end{equation} 
where $r:= (x_1^2 + \cdots + x_n^2)^{1/2}$, with
\begin{equation}
{\cal Y}^A_{|\eta| - 2m,\cdot}(x) =
{1 \over 4^m m! \Big ( n/2 + n(n-1)/2\alpha + |\eta| - 2m \Big )_m}
\tilde{T}_{|\eta|-2m}^A\Big ( {\Delta_A^m} E_\eta(x) \Big ),
\end{equation}
\begin{equation}
\tilde{T}_{k}^A = \sum_{j=0}^{[k/2]} {r^{2j} \over 4^j j! \Big (
-n/2 - n(n-1)/2\alpha - k + 2)_j} {\Delta_A^j}.
\end{equation}

Dunkl's theory also allows the non-symmetric Hermite polynomials to be
expressed in terms of the harmonic polynomials. Thus from
\cite[Prop.~3.9]{dunkl91a} we know that
\begin{equation}
e^{-{\Delta}_A/4} |r|^{2j} {\cal Y}^A_{k,\cdot}(x) = (-1)^j j!
L_j^{k + n(n-1)/2\alpha + n/2 - 1}(r^2) {\cal Y}^A_{k,\cdot}(x),
\end{equation}
where $L_j^a$ denotes the classical (one-variable) Laguerre polynomial,
and so applying the formula (\ref{exp.g}) to (\ref{ey}) we obtain
\begin{equation}\label{harmg}
E_\eta^{(H)}(x) = \sum_{m=0}^{[|\eta|/2]} (-1)^m m!
L_m^{(|\eta| - 2m +  n(n-1)/2\alpha + n/2 -1)} (r^2)
{\cal Y}^A_{|\eta| - 2m,\cdot}(x).
\end{equation}

\setcounter{equation}{0}
\section{The Laguerre Case}

In this section we investigate the Laguerre case, in an analogous
manner to the Hermite case of the previous section.
We begin by reviewing some
results in \cite{forr96c}. Firstly, we have
the type $B$ Dunkl operators
\begin{equation}
T_i^{(B)} := \dif{x_i} + {1 \over \alpha} \sum_{p \neq i}
\left( { 1 - s_{ip} \over x_i - x_p} + {1 - \sigma_i \sigma_p s_{ip} 
\over x_i + x_p}
\right) + {a + 1/2 \over x_i}(1 - \sigma_i)
\end{equation}
where $\sigma_j$ is the operator which replaces $x_j$ by $-x_j$. Instead
of working with the operator $T_i^{(B)}$ directly, in \cite{forr96c} we
found it convenient to work with the operator $B_i:=\frac{1}{4}
(T_i^{(B)})^2$ acting on functions of $x^2$, since in this case
$$
B_i = x_i^2\, \hat{T}_i^2 + (a+1) \hat{T}_i +
{1 \over \alpha} \sum_{p \ne i} s_{ip} \hat{T}_i
$$
where $\hat{T}_i$ is the type $A$ Dunkl operator in the variables $x^2$:
$$
\hat{T}_i = {1 \over 2x_i} \dif{x_i}
+ {1 \over \alpha} \sum_{p \ne i} {1 - s_{ip} \over x_i^2 - x_p^2}.
$$
Moreover, the operators $B_i$ enjoy the following property
\cite[Lemma 4.2]{forr96c}
\begin{lemma}\label{glow}
Let
\begin{equation}
\hat{\xi}_j := \alpha x_j^2\,\hat{T}_j + (1 - n) + \sum_{p > j} s_{jp}
\end{equation}
be the Cherednik operator (\ref{cherednik.1}) with the 
substitution $x_j \rightarrow x_j^2$ ($j=1,\dots,N$). Then we have
\begin{eqnarray*}
[ \hat{\xi}_j,B_i] & = & B_i s_{ij}, \quad i<j \\
{[} \hat{\xi}_j,B_i] & = & B_j s_{ij}, \quad i>j \\
{[} \hat{\xi}_j,B_j]& = &- \alpha B_j - \sum_{p < j} s_{jp} B_j - \sum_{p >
j} B_j s_{jp}.
\end{eqnarray*}
\end{lemma}

The non-symmetric Laguerre polynomials $\el_{\eta}(x^2)$, which are
eigenfunctions of (\ref{h-lag}) with eigenvalue $-4|\eta|$, are also
eigenfunctions of the operators \cite{forr96c}
\begin{equation}
l_i = \hat{\xi}_i - \al\,B_i
\end{equation}
which are mutually commuting and self-adjoint with respect to the
inner product (\ref{inn.l}). As in the Jack and Hermite cases, the
operators $l_i$, $s_j$ obey the relations
$$
l_i\,s_i - s_i \,l_{i+1} = 1, \qquad l_{i+1}\,s_i -s_i\, l_i = -1, \qquad
{}[l_i,s_j] = 0, \quad j\neq i, i+1
$$
which immediately gives
\begin{equation} \label{ammonia}
s_i\,\el_{\eta} = \left\{ \begin{array}{ll}
\bfrac{1}{\de_{i,\eta}}\,\el_{\eta} + \left( 1- \bfrac{1}{\de_{i,\eta}^2}
\right) \el_{s_i\eta} & \eta_i > \eta_{i+1} \\
\el_{\eta} & \eta_i = \eta_{i+1} \\
\bfrac{1}{\de_{i,\eta}}\,\el_{\eta} + \el_{s_i\eta} & \eta_i < \eta_{i+1}
\end{array} \right.
\end{equation}

Let $\Psi$ be the analogue of the operator $\Phi$ acting on functions of
$x^2$. That is, $\Psi:=x_n^2\,s_{n-1}\cdots s_2\,s_1$ which clearly acts
on non-symmetric Jack functions by
$$
\Psi\,E_{\eta}(x^2) = E_{\Phi\eta}(x^2)
$$
where, as before $\Phi\eta = (\eta_2,\eta_3,\ldots,\eta_n,\eta_1+1)$.

We also introduce the non-symmetric analogue of the generalized
factorial function as
\begin{eqnarray}
[c]_{\eta} &:=& \prod_{s\in\eta}\left( c + a'(s) -l'(s)/\al
\right)\nonumber\\
&=& \prod_{i=1}^n \left( c -\eta_i +\bar{\eta}_i/\al \right)_{\eta_i^+}
\end{eqnarray}
where $a'(s)$, $l'(s)$ are defined in (\ref{guion}). Following the
arguments of \cite[Lemmas 4.1, 4.2]{sahi96a} (indeed
$e_{\eta} = \al^{|\eta|}\: [n/\al +1]_{\eta}$) the coefficients
$[c]_{\eta}$ have the properties
\begin{equation} \label{codo}
[c]_{s_i\eta} = [c]_{\eta} \quad \mbox{for all $i$}, \qquad
\frac{[c]_{\Phi\eta}}{[c]_{\eta}} = c + \bar{\eta}_1/\al, \qquad
\frac{[c]_{\eta}}{[c]_{\ph\eta}} = c -1 + \bar{\eta}_n/\al .
\end{equation}
{}From the first property we see that
\begin{equation}\label{fact}
[c]_\eta = [c]_{\eta^+} = 
\prod_{i=1}^n \Big ( c - {i-1 \over \alpha} \Big )_{\eta_i^+}.
\end{equation}

Finally let $\widehat{\Psi} := B_1\,s_1\,s_2\cdots s_{n-1}$. Then
Lemma \ref{unfettered} and Corollory \ref{antemano} have the following
analogue in the Laguerre case

\begin{lemma}\label{lemq}
We have
\begin{eqnarray*}
(a)\hspace{3cm} \hat{\xi}_j\,\widehat{\Psi} &=& \widehat{\Psi}\,
\hat{\xi}_{j-1} \hspace{3cm} 2\leq j \leq n \\
(b)\hspace{3cm} \hat{\xi}_1\,\widehat{\Psi} &=& \widehat{\Psi}\,
\left(\hat{\xi}_{n} -\al \right) \\
(c)\hspace{2.1cm}\widehat{\Psi} \, E_{\eta}(x^2) &=& 
\frac{1}{\al}\;\frac{[a+q]_{\eta}}
{[a+q]_{\ph\eta}}\;\frac{d'_{\eta}}{d'_{\ph\eta}}\;
E_{\ph\eta}(x^2)
\end{eqnarray*}
where $q:=1 +(n-1)/\al$.
\end{lemma}
{\it Proof.}\quad (a) and (b) follow immediately from Lemma
\ref{glow}, in analogy to the corresponding results in the
previous section. To prove (c), we know from (a) and (b) that 
$$
\widehat{\Psi} \, E_{\eta}(x^2) = c_{\eta}\;
E_{\ph\eta}(x^2)
$$
for some constant $c_{\eta}$.
An examination of the leading term shows that
$$
c_{\eta} = \left( a+q-1 +\bar{\eta}_n/\al \right)
\left(\frac{\bar{\eta}_n +n-1}{\al} \right)
$$
This can be simplified with (\ref{jude}) and (\ref{codo})
to give the desired result. \hfill $\Box$

\vspace{2mm}
We now construct raising and lowering operators for the non-symmetric
Laguerre polynomials in the same manner as for the Hermite case.
{}From Dunkl \cite{dunkl91a} we have the relations
\begin{equation} \label{membrete}
[x_i,\Delta_B] = -2\,T_i^{(B)}, \hspace{3cm}
(T_i^{(B)})^{\star} = 2x_i - T_i^{(B)} ,
\end{equation}
where $\Delta_B := \sum_i (T_i^{(B)})^2$ and
$(T_i^{(B)})^{\star}$ denotes the adjoint of $T_i^{(B)}$ with
respect to the Laguerre inner product (\ref{inn.l}), while from
\cite{forr96c} we have
$$
[\hat{\xi}_i, \Delta_B ] = -4\al B_i 
$$
Using (\ref{membrete}) we have the following expression for $\shs$:
\begin{eqnarray}
\shs = s_{n-1}\cdots s_1\; B_1^{\star} &=& \frac{1}{4}
s_{n-1}\cdots s_1\; \left(2x_1 - T_1^{(B)} \right)^2 \nonumber\\
&=& \Psi + \sfrac{1}{4}[\Psi,\Delta_B] +
s_{n-1}\cdots s_1\; B_1 \label{disparates}
\end{eqnarray}
The Laguerre analogue of Proposition \ref{acotar} is given by
\begin{prop}
The operators $l_i$ satisfy
\begin{eqnarray*}
(a)\qquad l_n\;\shs &=& \shs\;(l_1+\al) \hspace{2cm}
l_i\;\shs = \shs\;l_{i+1} \qquad 1\leq i\leq n-1 \qquad\\
(b)\qquad l_1\;\widehat{\Psi}\;&=& \widehat{\Psi}\;(l_n-\al) \hspace{2cm}
l_i\;\widehat{\Psi} = \widehat{\Psi}\;l_{i-1} \qquad 2\leq i\leq n
\end{eqnarray*}
\end{prop}
{\it Proof.}\quad We prove only the first equation appearing in (a), the
others being similar. Using the representation (\ref{disparates}) we have
\begin{eqnarray*}
B_i\,\shs &=& \shs\,B_{i+1} +2\tau -[\Delta_B,\tau] \\
\hat{\xi}_i\,\shs &=& \hat{\xi}_{i+1}\,\shs +2\al\tau +
s_{in}\,s_{n-1}\cdots s_1\,B_1 + s_{n-1}\cdots s_{i+1}\,s_{i-1}\cdots
s_1\,B_{i+1}
\end{eqnarray*}
where
\begin{eqnarray*}
\tau  &:=& \sfrac{1}{4}\left(B_i\,\Psi - \Psi\, B_{i+1} \right) 
= \sfrac{1}{4}s_{n-1}\cdots s_1\,[B_{i+1},x_1^2] \\
&=& \sfrac{1}{4}s_{n-1}\cdots s_1\,\left( x_{i+1}^2\hat{T}_{i+1}
+ x_1^2\hat{T}_1 + a+1 +\sum_{p\neq i+1}s_{ip} \right)\left(
-\frac{1}{\al}\,s_{i+1,1} \right)
\end{eqnarray*}
It remains to calculate the commutator $[\Delta_B,\tau]$. To this end
we recall that when acting on functions of $x^2$, we have
$T_i^{(B)} = 2x_i\hat{T}_i$ whence from (\ref{membrete})
$$
[\Delta_B,x_i^2\,\hat{T}_i] = \sfrac{1}{2}[\Delta_B,x_i\,T_i^{(B)}]
=4\,B_i\;.
$$
Thus
$$
[\Delta_B,\tau] = -\frac{1}{\al}\,s_{n-1}\cdots s_1\,
s_{i+1,1}\,(B_{i+1}+B_1)
$$
and the result now follows using the fact that
$$
s_{n-1}\cdots s_1\,s_{i+1,1} = s_{n-1}\cdots s_{i+1}\,s_{i-1}\cdots
s_1 = s_{in}\,s_{n-1}\cdots s_1
$$
\hfill $\Box$

\begin{cor}
The operators $\widehat{\Psi}$ and $\shs$ act on the non-symmetric
Laguerre polynomials via
\begin{eqnarray*}
\widehat{\Psi}\,\el_{\eta}(x^2) &=& \frac{1}{\al}\;\frac{[a+q]_{\eta}}
{[a+q]_{\ph\eta}}\;\frac{d'_{\eta}}{d'_{\ph\eta}}\;\el_{\ph\eta}(x^2) \\
\shs\,\el_{\eta}(x^2) &=& \el_{\Phi\eta}(x^2)
\end{eqnarray*}
\end{cor}

We can now use the above result to calculate the norm of the functions
$\el_{\eta}(x^2)$

\begin{prop}\label{norm.l}
We have
$$
\innl{\el_{\eta}}{\el_{\eta}} = \frac{[a+q]_{\eta}}{\al^{|\eta|}}\,
\frac{d'_{\eta}e_{\eta}}{d_{\eta}}\;{\cal N}^{(L)}_0
$$
where
$$
{\cal N}^{(L)}_0 := \innl{1}{1} =
\alpha^{(1-N-(N-1)^2/\alpha)} \prod_{j=0}^{N-1}
{\Gamma(1+(j+1)/\alpha) \Gamma(a+1+j/\alpha) \over \Gamma(1+1/\alpha)}.
$$
\end{prop}
{\it Proof.}\quad Similar to the proof of Proposition \ref{tierra}.
\hfill $\Box$

\subsection{Generating function}

Let
\begin{equation}\label{kb}
{\cal K}_B(x^2;y^2) = \sum_{\eta} \frac{\al^{|\eta|}}{[a+q]_{\eta}}\;
\frac{d_{\eta}}{d'_{\eta}\,e_{\eta}} \;E_{\eta}(x^2)\,E_{\eta}(y^2)
\end{equation}
where $q$ is defined as in Lemma \ref{lemq}.
Then the following result is proved in the manner of Theorem \ref{main}

\begin{thm} \label{mainB}
The function ${\cal K}_B(x^2;y^2)$ possesses the following properties
\begin{eqnarray*}
(a)\qquad s_i^{(y)}\;{\cal K}_B(x^2;y^2) &=& s_i^{(x)}\;{\cal K}_B(x^2;y^2) \\
(b)\qquad \widehat{\Psi}^{(y)}\;{\cal K}_B(x^2;y^2) &=& 
\Psi^{(x)}\;{\cal K}_B(x^2;y^2) \\
(c)\qquad B_i^{(y)}\;{\cal K}_B(x^2;y^2) &=& x^2_i\;{\cal K}_B(x^2;y^2) 
\end{eqnarray*}
\end{thm}
The formula (c) of Theorem \ref{mainB} implies that
\begin{equation}\label{DelB}
\Delta_B \, {\cal K}_B(x^2;y^2) = 4 p_1(x^2)  {\cal K}_B(x^2;y^2).
\end{equation}
{}From (\ref{DelB}) and the fact that $H^{(L)} = \Delta_B -2 \tilde{E}_1$,
we can proceed as in the proof of Proposition \ref{hermgf} to
derive a generating function for the non-symmetric Laguerre polynomials.

\begin{prop}\label{laggf}
We have
\begin{equation}\label{genl}
\sum_\eta {(-\alpha)^{|\eta|} \over [a+q]_\eta} {d_\eta \over d_\eta' e_\eta}
E_\eta^{(L)}(x) E_\eta(z) = {\cal K}_B(x;-z) e^{p_1(z)}.
\end{equation}
\end{prop}

Analogous to the case of the symmetric Laguerre polynomials 
\cite[Prop.~4.3]{forr96a}
we can
use the generating function (\ref{genl}) to express the non-symmetric
Laguerre polynomials as a series in non-symmetric Jack polynomials
involving the generalized binomial coefficients.
\begin{prop}
We have
\begin{equation}\label{f1}
E_\eta^{(L)}(x) = {(-1)^{|\eta|}[a+q]_\eta e_\eta \over d_\eta} 
\sum_{\nu} {(-1)^{|\nu|} \over [a+q]_\nu}{d_\nu \over e_\nu}
\left ( { \eta \atop \nu} \right ) E_\nu(x)
\end{equation}
\begin{equation}\label{f2}
E_\eta(x) = {[a+q]_\eta e_\eta \over d_\eta} \sum_{\nu}
 {1\over [a+q]_\nu}{d_\nu \over e_\nu}
 \left ( { \eta \atop \nu} \right ) E_\nu^{(L)}(x)
\end{equation}
\end{prop}

\vspace{.2cm}
\noindent
{\it Proof.}\quad
The formula (\ref{f1}) follows from (\ref{genl}) by applying the identity
(\ref{n.conn}) on the right hand side and equating coefficients of
$E_\eta(z)$. The formula (\ref{f2}) follows from  (\ref{genl}) by
multiplying both sides by $e^{-p_1(z)}$, using the identity (\ref{n.conn})
with $z$ replaced by $-z$ on the left hand side and equating like
coefficients of $E_\eta(z)$.
\vspace{.2cm}

By substituting $x=0$ in (\ref{f1}) we obtain the explicit value of
$E_\eta^{(L)}$ at the origin.
\begin{cor}\label{cor9}
We have
$$
E_\eta^{(L)}(0) = {(-1)^{|\eta|} [a+q]_\eta e_\eta \over d_\eta}.
$$
\end{cor}

\subsection{Other generating functions}
For the symmetric Laguerre polynomials we presented \cite{forr96a}
three generating functions including the analogue of (\ref{genl}). These
generating functions all rely on a partial differential equation
\cite[Prop.~A.1]{forr96a} for the generalized hypergeometric
function
$$
{}_2{\cal F}_1(x;y) := \sum_{\kappa}{1 \over |\kappa|!}
{[a]_\kappa [b]_\kappa \over [c]_\kappa} {C_\kappa(x) C_\kappa(y)
\over C_\kappa(1^n)}.
$$
A similar partial differential equation is satisfied by the function
\begin{equation}\label{2K1}
{}_2{\cal K}_1(x;y) := \sum_{\eta} \alpha^{|\eta|}
{[a]_\eta [b]_\eta \over [c]_\eta}{d_\eta \over d_\eta' e_\eta}
E_\eta(x) E_\eta(y).
\end{equation}

\begin{prop}\label{pde1}
The function ${}_2{\cal K}_1(x;y)$ satisfies the p.d.e.
\begin{equation}\label{pde}
\tilde{D}_1^{(x)}F + \Big ( c - {N-1 \over \alpha} \Big )
\tilde{E}_0^{(x)}F -  \Big ( a+b  - {N-1 \over \alpha} \Big )
\tilde{E}_2^{(y)}F - {1 \over 2}[\tilde{D}_2,\tilde{E}_2]^{(y)}F =
abp_1(y)F
\end{equation}
where $\tilde{E}_0$, $\tilde{E}_2$ are defined by (\ref{ep}),
 $\tilde{D}_1$ is defined by (\ref{dt1}) and  $\tilde{D}_2$ by
  (\ref{dt2}). In fact ${}_2{\cal K}_1$ is the unique solution
of this p.d.e.~of the form
\begin{equation}\label{trial}
F(x,y) =  \sum_{\eta} A_{\eta^+} {d_\eta \over d_\eta' e_\eta}
E_\eta(x) E_\eta(y), \qquad (A_0 = 1)
\end{equation}
\end{prop}

\vspace{.2cm}
\noindent
{\it Proof.} \quad Substituting (\ref{trial}) in (\ref{pde}), we see
that the action of all the operators on $E_\eta(x)$ and $ E_\eta(y)$
is specified in Proposition \ref{action.b}, while $p_1(y) E_\eta(y)$ can be
written according to (\ref{pex}). Equating like coefficients of
$E_\eta(x)E_\eta(y)$, $|\nu| = |\eta| +1$ gives
\begin{eqnarray}\lefteqn{
\Big ( {\epsilon_\nu - \epsilon_\eta \over 2} + c -{N-1 \over \alpha}
\Big ) \left ( {\nu \atop \eta} \right ) A_{\nu^+}} \nonumber \\
&& =\bigg ( \Big (  {\epsilon_\nu - \epsilon_\eta \over 2} -
{N - 1 \over \alpha} \Big )\Big ( a+b - {N - 1 \over \alpha}
+  {\epsilon_\nu - \epsilon_\eta \over 2} \Big ) + ab \bigg )
 \left ( {\nu \atop \eta} \right ) A_{\eta^+} \label{abov}
 \end{eqnarray}
Let's suppose that $\eta_1 \le \eta_2  \cdots \le \eta_N$ and
$\nu_{N-i+1} = \eta_{N-i+1}+1$, $\nu_j = \eta_j$ $j=1,\dots,n \:
(j \ne i)$. Suppose
furthermore that $ \eta_{N-i+1}+1 \le  \eta_{N-i+2}$. From the ordering
in the expansion (\ref{1.3}) and the fact that coefficients $a_{\eta
\nu}$ therein are positive, we see from the 
definition (\ref{binomial}) that $ \left ( {\nu \atop \eta} \right )
$ is non-zero, and so can be cancelled from (\ref{abov}). Noting
$\eta_j^+ = \eta_{N-j+1}$, we see from (\ref{He}) that
$ {\epsilon_\nu - \epsilon_\eta \over 2}=\eta_i^+ + {N-i \over \alpha}$,
and so
$$
\Big ( c + \eta_i^+ - {i-1 \over \alpha} \Big ) A_{\nu^+} =
\alpha \Big (  a + \eta_i^+ - {i-1 \over \alpha} \Big )
 \Big (  b + \eta_i^+ - {i-1 \over \alpha} \Big ) A_{\eta^+}.
$$
This is a first order difference equation in the parts of the
partition $\eta^+$ and so, for a given initial
condition ($A_0 = 1$), has a unique solution. It is straightforward to
verify from (\ref{fact}) that the solution is $A_{\eta^+} =
[a]_\eta[b]_\eta / [c]_\eta$, as required. \hfill$\Box$

\vspace{.2cm}
Of particular interest is  Proposition \ref{pde1} with the change
of variables $y \mapsto y/b$ and $b \to \infty$. This gives
\begin{cor}\label{c5}
The p.d.e.
\begin{equation}\label{pde2}
\tilde{D}_1^{(x)} F + \Big ( c - {N-1 \over \alpha} \Big )
\tilde{E}_0^{(x)} F - \tilde{E}_2^{(y)} F = a p_1(y) F
\end{equation}
is satisfied by
$$
F = {}_1{\cal K}_1(a;c;x;y)
$$
where
$$
 {}_1{\cal K}_1(a;c;x;y) := \sum_\eta \alpha^{|\eta|}{[a]_\eta \over [c]_\eta}
 {d_\eta \over d_\eta' e_\eta} E_\eta(x) E_\eta(y).
$$
\end{cor}

By proceeding as in the derivation of \cite[Prop.~4.2]{forr96a} we
can deduce from Corollary \ref{c5} a generating function for the
$E_\eta^{(L)}$ involving the function ${}_1{\cal K}_1$.

\begin{prop}\label{p5}
We have
\begin{equation}\label{p6}
\prod_{i=1}^n (1 - z_i)^{-c-q} {}_1 {\cal K}_1(c+q;a+q;-x;{z \over 1
-z}) = \sum_{\eta} (-\alpha)^{|\eta|}{[c+q]_\eta \over [a+q]_\eta} {d_\eta \over 
d_\eta' e_\eta}  E_\eta^{(L)}(x) E_\eta(z).
\end{equation}
\end{prop}

Replacing $z$ by $z/c$ in (\ref{p6}) 
and taking the limit $c \to \infty$ we reclaim
the generating function (\ref{genl}), while setting $a=c$ the following
generating function results.

\begin{prop}
We have
\begin{equation}\label{g4}
\prod_{i=1}^n (1 - z_i)^{-a-q} {\cal K}_A(-x;{z \over 1
-z}) = \sum_{\eta} (-\alpha)^{|\eta|} {d_\eta \over
d_\eta' e_\eta} E_\eta^{(L)}(x) E_\eta(z),
\end{equation}
where ${\cal K}_A$ is defined by (\ref{ker-A}).
\end{prop}
As in the case of the symmetric Laguerre polynomials, by using the
generating function (\ref{g4}), the orthogonality of $\{E_\eta^{(L)}\}$
with respect to (\ref{inn.l}) and Corollary \ref{cor9} we can provide
an alternative proof of Proposition \ref{norm.l} (see \cite[Prop.~4.1]
{forr96a} for details).

\subsection{Generalized Hankel transform}
The generating function formulas (\ref{genl}), (\ref{p6}) and (\ref{g4})
are all direct analogues of the generating functions for the symmetric
Laguerre polynomials. In the symmetric case, the analogue of the
function ${\cal K}_B(x;y)$ also occurs as the kernel in an integral transform
which relates the symmetric Laguerre and Jack polynomials
\cite[Prop.~4.11 and Cor.~4.1]{forr96a}. The non-symmetric analogue is
easily obtained by following a similar strategy.

\begin{prop}
Let
\begin{equation}
d\mu^{(L)}(y) := \prod_{j=1}^n y_j^a e^{-y_j} 
\prod_{1 \le j < k \le n} |y_k - y_j|^{2/\alpha} dy_1 \dots dy_n.
\end{equation}
We have
\begin{equation}
\int_{[0,\infty)^n} {\cal K}_B(x;-z_a)  {\cal K}_B(x;-z_b)d\mu^{(L)}(x)
= {\cal N}_0^{(L)} e^{-p_1(z_a)} e^{-p_1(z_b)}  {\cal K}_B(z_a;z_b).
\end{equation}
\begin{equation}\label{int.ll}
\int_{[0,\infty)^n} {\cal K}_B(x;-z_a) E_\eta^{(L)}(x) d\mu^{(L)}(x)
=  {\cal N}_0^{(L)} e^{-p_1(z_a)} E_\eta(-z_a).
\end{equation}
\begin{equation}\label{int.l}
\int_{[0,\infty)^n} {\cal K}_B(x;-z_a) E_\eta(-x) d\mu^{(L)}(x)
=  {\cal N}_0^{(L)} e^{-p_1(z_a)} E_\eta^{(L)}(z_a).
\end{equation}
\end{prop}

Using (\ref{int.l}) it is straightforward to derive the Laguerre analogue
of the summation formula in Proposition \ref{green.h}
(see \cite[Prop.~4.12]{forr96a} for the symmetric
case and further details).

\begin{prop}\label{aft}
For $|t| < 1$ we have
\begin{eqnarray*}
\sum_{\eta} \frac{1}{{\cal N}_{\eta}^{(L)}}\;
E_{\eta}^{(L)}(x)\:E_{\eta}^{(L)}(y) t^{|\eta|}
= \frac{1}{{\cal N}_0^{(L)}} (1-t)^{-N(a+q)} \hspace{4cm}\\
\times\exp\left(-\frac{t}{1-t}\left(p_1(x) + p_1(y)\right)\right)\;
{\cal K}_B\left(
{y \over 1-t}; {tx \over 1 - t} \right )\end{eqnarray*}
\end{prop}
We can use Proposition \ref{aft} to prove the analogue of
the asymptotic expansion (\ref{3.30}).
\begin{prop}\label{asyb}
We have
$$
{\cal K}_{B}\Big ({x^2 \over 2 \tau};{y^2 \over 2 \tau} \Big )
\: \sim \: {\pi^{-n/2} 2^{n(a+1/2) + n(n-1)/\alpha} {\cal N}_0^{(L)}
\over ( \prod_{1 \le j < k \le n} (x_j^2 - x_k^2) (y_j^2 - y_k^2)/
\tau)^{1/\alpha}} \prod_{j=1}^n (x_j y_j / \tau)^{-(a+1/2)}
e^{x_j y_j /\tau}
$$
as $\tau \to 0$.
\end{prop}

\vspace{.2cm}
\noindent
{\it Proof.} \quad This follows from the interpretation of the sum in
Proposition \ref{aft} as the Green function
for a Schr\"odinger equation. See the derivation of the
symmetric counterpart of this result \cite[eq.~(5.47)]{forr96a} for
further details (in fact, as is the case in the type $A$ case (\ref{3.30}),
the asymptotic expansion is identical with its symmetric counterpart).

\vspace{.2cm}
To further develop ${\cal K}_{B}$ as an integral kernel we will require
the exponential operator formula \linebreak
\cite[eqs.~(4.2)\& (4.4)]{forr96c}
\begin{equation}\label{exp.l}
E_\eta^{(L)}(x^2) = e^{-\Delta_B/4} E_\eta(x^2),
\end{equation}
which can also be derived from (\ref{DelB}) and the generating function
(\ref{genl}). Substituting (\ref{exp.l}) in (\ref{int.l}) and
(\ref{int.ll}) with the change of variables $x \mapsto x^2$,
$z_a \mapsto z^2$, and using the fact that $\{E_\eta^{(L)}\}$ is a basis
for analytic functions, we obtain the following formulas.
\begin{prop}
We have
\begin{equation}\label{eq.a}
\int_{-(\infty,\infty)^n} {\cal K}_B(x^2;-z^2) \Big (
 e^{-\Delta_B/4} f(x^2) \Big ) \, d\mu^{(L)}(x^2) = 
 {\cal N}_0^{(L)} e^{-p_2(z)} f(-z^2)
 \end{equation}
\begin{equation}\label{eq.b}
\int_{-(\infty,\infty)^n} {\cal K}_B(x^2;-z^2)f(-x^2) \,  d\mu^{(L)}(x^2) =
 {\cal N}_0^{(L)} e^{-p_2(z)}  e^{-\Delta_B/4} f(z^2)
 \end{equation}
 where $f$ is an analytic function such that all terms converge.
 \end{prop}

 Analogous to the pairing (\ref{pairing}) we define the $B$-type pairing
 \begin{equation}\label{pair.l}
 [p,q]_L = p((T^{(B)})^x) q(x),
 \end{equation}
 where $p$ and $q$ are homogeneous polynomials of degree $|\kappa|$ say.
 According to the theory of Dunkl \cite[Th.~3.10]{dunkl91a},
(\ref{pair.l}) is related to the exponential operator in 
(\ref{exp.l}) by
\begin{equation}
[p(x^2),q(x^2)]_L = {2^{2|\kappa|}
\over {\cal N}_0^{(L)}} \int_{(-\infty,\infty)^n}
\Big (  e^{-\Delta_B/4} p(x^2) \Big )
\Big (  e^{-\Delta_B/4} q(x^2) \Big ) \,  d\mu^{(L)}(x^2).
\end{equation}
Using (\ref{exp.l}), the orthogonality of $\{E_\eta^{(L)}\}$ with
respect to (\ref{inn.l}) and Proposition \ref{norm.l} we obtain

\begin{prop}
We have
\begin{equation}\label{pairingb}
[E_\nu(x^2),E_\eta(x^2)]_L = 2^{2|\kappa|}{[a+q]_\eta \over \alpha^{|\eta|}} 
{d_\eta' e_\eta \over d_\eta} \delta_{\nu,\eta}.
\end{equation}
\end{prop}

Another derivation of (\ref{pairingb}) can be deduced from a recent
work of Dunkl \cite{dunkl96b} on intertwining operators of type $B$.
The key formula \cite[proof of Th.~4.2]{dunkl96b} is that for
$f \in A_\lambda$ (this is the space defined above (\ref{vf}))
\begin{equation} \label{vb}
V^{(B)} \xi f(x^2) = {1 \over 2^{2|\lambda|} [n/\alpha +
1]_\lambda^{(\alpha)} [a+q]_\lambda^{(\alpha)}}\;f(x^2)
\end{equation}
with $\xi$ defined below (\ref{vf}) and $V^{(B)}$ defined by the 
intertwining relation $T_i^{(B)} V^{(B)} = V^{(B)} { \partial \over
\partial x_i}$ and the normalization $V \, 1 = 1$. By following the
working which led to (\ref{pfga}) we deduce from (\ref{vb}) that
for $f,g \in A_\lambda$
\begin{equation}
[f(x^2), g(x^2)]_L = 2^{2|\lambda|} [n/\alpha + 1]_\lambda^{(\alpha)}
[a+q]_\lambda^{(\alpha)}  \sum_{\nu'} f_{\nu'} \sum_\nu A_{\nu'
  \nu} g_\nu.
\end{equation}
Substituting (\ref{sip}) and using the facts noted below (\ref{dufgh})
reclaims (\ref{pairingb}).

 The formula (\ref{eq.a}) has occured in the type-$B$ theory of Dunkl
 \cite[Prop.~2.1]{dunkl92a} with ${\cal K}_B$ replaced by a certain 
 kernel $K_B(
 \sqrt{2}x,\sqrt{2}z)$, $f(y)$ by $f(\sqrt{2}y)$, $x$ by $x/\sqrt{2}$
 and $z$ by $z/\sqrt{2}$. Thus we must have
 \begin{equation}\label{e9}
 K_B(x,z) = {\cal K}_B(x^2/2;z^2/2).
 \end{equation}
 We remark also that Dunkl \cite{dunkl91a} has proved, without using the
 explicit formula (\ref{e9}), that $K_B$ has the properties (a) and
 (b) of Theorem \ref{mainB}. Furthermore, it follows from the properties of $K_B$
 established in  \cite{dunkl91a} that $ {\cal K}_B$ is an entire
 function of all variables, and satisfies the uniform bound
  \begin{equation}
| {\cal K}_B(x^2;z^2)| \le e^{cp_1(|x|) p_1(|y|)},
\end{equation}
for some $c > 0$ (c.f.~Proposition \ref{asyb}).

\subsection{The generalized Hankel transform and its relation to the
Laplace transform}
In refs.~\cite{macunp1,yan92a} (see also 
\cite{forr96a}) the generalized Hankel transform ${\cal H}$ is defined
as the symmetric version of
 \begin{equation}\label{hank}
 {\cal  H}[f(x^2)](z^2) = {1 \over {\cal N}_0^{(L)}} \int_{(-\infty, \infty)^n}
   {\cal K}_B(x^2;-z^2) f(x^2) \, d\mu^{(L)}(x^2),
    \end{equation}
(i.e.~(\ref{hank}) with ${\cal K}_B(x^2;-z^2)$ replaced by ${}_0{\cal F}_1
(a+q;x^2;-z^2)$).
In ref.~\cite{macunp1} it is shown, on the basis of some conjectures
concerning the generalized Laplace transform, that the symmetric Hankel
transform is an isometry with respect to the inner product (\ref{inn.h}),
and is further related to the Laplace transform by a generalization of
Tricomi's theorem. In ref.~\cite{yan92a} these properties are proved in
the case $n=2$.

In the non-symmetric case Dunkl \cite{dunkl91a,dunkl92a} has defined the
Hankel transform by (\ref{hank}) with $z^2$ replaced by $iz^2$. The 
isometry property is proved, and the relationship between  the
Hankel transform and polynomials annihilated by the operator
$\Delta_{B}$ explored in some detail. Here we complement Dunkl's
theory by extending the results
of ref.~\cite{macunp1,yan92a} (see also \cite{forr96a}) for the symmetric case.

First we will calculate the generalized Laplace transform
(\ref{laplace}) of $\prod_{j=1}^n x_j^a \, E_\eta^{(L)}$ and use the
result to calculate the Laplace transform of $\prod_{j=1}^n
 x_j^a \, E_\eta$.

 \begin{prop}
We have
\begin{equation}\label{l.el}
{\cal L}\Big ( \prod_{j=1}^n x_j^a \, E_\eta^{(L)}(x) \Big ) =
[a+q]_\eta \, {\cal N}_0^{(L)}  \prod_{j=1}^n t_j^{-(a+q)}
E_\eta(1 - {1 \over t})
\end{equation}
\begin{equation}\label{l.e}
{\cal L}\Big ( \prod_{j=1}^n x_j^a \, E_\eta(x) \Big ) =
[a+q]_\eta \, {\cal N}_0^{(L)} \prod_{j=1}^n t_j^{-(a+q)}
E_\eta^{(L)}( {1 \over t}).
\end{equation}
\end{prop}

\vspace{.2cm}
\noindent
{\it Proof.}\quad 
To derive (\ref{l.el}) we multiply both sides of the generating function
(\ref{g4}) by $E_\eta^{(L)}(x)$ and integrate with respect to the
measure $d\mu^{(L)}(x)$ over the region $[0,\infty)^n$. This gives
\begin{equation}\label{sun}
\prod_{j=1}^n (1 - z_j)^{-(a+q)} \int_{[0,\infty)^n}
{\cal K}_A(-x;{z \over 1 - z}) E_\eta^{(L)}(x) \, d\mu^{(L)}(x)
= {\cal N}_0^{(L)} [a+q]_\eta 
 E_\eta(z)
\end{equation}
where on the r.h.s.~we have used the
orthogonality of $\{E_\eta^{(L)} \}$ with respect
to (\ref{inn.l}) and the normalization in Proposition \ref{norm.l}.
Using the identity (\ref{pK}) we see that
$$
{\cal K}_A(-x;{z \over 1 - z})  \, d\mu^{(L)}(x)
= \prod_{j=1}^n x_j^a {\cal K}_A(-x;{1 \over 1 - z}) dx_1 \dots dx_n.
$$
Thus, with $1/(1-z) = t$,
the integral on the left hand side of (\ref{sun}) is the generalized
Laplace transform in (\ref{l.el}), and the first result follows.
The second result follows from the first by replacing $\eta$ by
$\nu$, and summing over $\nu$ with appropriate $\nu$ dependent factors
given in (\ref{f1}) so that $E_\eta^{(L)}(x)$ on the left hand side becomes
$ E_\eta(x)$. Performing the same operation on the r.h.s.~we see
from (\ref{binomial}) and (\ref{f1}) that the resulting sum is equal 
to $E^{(L)}_\eta(1/t)$,
as required. \hfill$\Box$

\vspace{.2cm}
Using (\ref{l.el}) we can compute the generalized Hankel transform
of ${\cal K}_A(-y^2;x^2)$.

\begin{prop}\label{kab}
We have 
\begin{equation}
{\cal H}[{\cal K}_A(-y^2;x^2)](z^2) = {\cal K}_A(-{1 \over y^2};z^2)
\prod_{j=1}^n y_j^{-2(a+q)},
\end{equation}
where $y$ is regarded as a parameter.
\end{prop}

\vspace{.2cm}
\noindent
{\it Proof.} \quad
{}From the definitions
$$
{\cal H}[{\cal K}_A(-y^2;x^2)](z^2) = {1 \over {\cal N}_0^{(L)}}
{\cal L}[\prod_{j=1}^n x_j^a \, {\cal K}_B(x;-z^2)](y^2).
$$
Computation of the right hand side using (\ref{l.el})
term-by-term on the series formula
(\ref{kb}) gives the stated result. \hfill$\Box$

\vspace{.2cm}
\noindent
Note from (\ref{gFt}) and (\ref{gFt.i}) that we can use
${\cal K}_A(-y^2;x^2)$ to form an arbitrary function of $x^2$
for which the generalized Fourier transform theory is applicable.
{}But from Proposition \ref{kab} we immediately have that
\begin{equation}
{\cal H}^2[{\cal K}_A(-y^2;x^2)](z^2) = {\cal K}_A(-y^2;z^2)
\end{equation}
which when combined with the above remark says that in general ${\cal H}$ is a
projection operator: ${\cal H}^2=1$. We also have that ${\cal H}$ is 
an isometry with respect to the inner product (\ref{inn.l}).
To see this, note that
\begin{eqnarray}\label{mon}
\langle {\cal K}_A(-y_1;x^2), {\cal K}_A(-y_2;x^2) \rangle_L & = &
\int_{(-\infty,\infty)^n} {\cal K}_A(-y_1;x^2) {\cal K}_A(-y_2;x^2) \,
d\mu^{(L)}(x^2) \nonumber \\
& = & \prod_{j=1}^n (y_1)_j^{-(a+q)} {}_1 {\cal K}_0(a+q;
- {1 \over y_1}; y_2) \nonumber \\
& = &  \prod_{j=1}^n (y_2)_j^{-(a+q)} {}_1 {\cal K}_0(a+q;
- {1 \over y_2}; y_1) 
\end{eqnarray}
where to obtain the second line we have used the series expansion 
(\ref{ker-A}) to replace ${\cal K}_A(-y_2;x^2)$
and integrated term-by-term using (\ref{l.e}), while to
obtain the final line the symmetry with respect to the interchange of
$y_1$ and $y_2$ has been used ($ {}_1 {\cal K}_0$ is given by
(\ref{2K1}) with $b=c$). On the other hand, from Proposition
\ref{kab} we have
\begin{equation}\label{fairway}
\langle {\cal H}[ {\cal K}_A(-y_1;x^2)],{\cal H}
[ {\cal K}_A(-y_2;x^2)] \rangle_L \!= \!
\prod_{j=1}^n (y_1)^{-(a+q)}_j(y_2)^{-(a+q)}_j\langle 
{\cal K}_A(-\frac{1}{y_1};x^2), {\cal K}_A(-\frac{1}{y_2};x^2) \rangle_L.
\end{equation}
Substituting the second last expression of (\ref{mon}) with $y_1$, $y_2$
replaced by $1/y_1$,  $1/y_2$ in the right hand side we see that the final
expression of (\ref{mon}) results, and so the right hand side 
of (\ref{fairway}) is precisely
$\langle {\cal K}_A(-y_1;x^2),{\cal K}_A(-y_2;x^2) \rangle_L$, as required.

\vspace{.2cm}

Let us now relate the generalized Laplace transform and Hankel transforms,
by giving the generalization of Tricomi's theorem.
\begin{prop}
The formula
\begin{equation}\label{mon1}
g(z) = {\cal H}[f(x^2)](z)
\end{equation}
holds if and only if
\begin{equation}\label{mon2}
{\cal L}[\prod_{j=1}^n z_j^a \, g(z)](t) = \prod_{j=1}^n t_j^{-(a+q)}
{\cal L}[\prod_{j=1}^n z_j^a \, f(z)]({1 \over t})
\end{equation}
\end{prop}

\vspace{.2cm}
\noindent
{\it Proof. }\quad First assume (\ref{mon1}). Since ${\cal H}^2 = 1$ we have
${\cal H}[g(x^2)](z) = f(z)$.
Use of this formula, the fact that ${\cal H}$ is an
isometry with respect to (\ref{inn.l}), and Proposition \ref{kab} gives
\begin{eqnarray*}
{\cal L}[\prod_{j=1}^n  z_j^a \, g(z)](t) & = &
\langle {\cal K}_A(-t;z), g(z) \rangle_L \\
& = & \langle {\cal H}[{\cal K}_A(-t;x^2)](z), {\cal H}[g(x^2)](z) \rangle_L \\
& = & \prod_{j=1}^n t_j^{-(a+q)} \langle
{\cal K}_A(-{1 \over t}; z ), f(z) \rangle_L,
\end{eqnarray*}
which is equivalent to (\ref{mon2}). Now assume (\ref{mon2}). Proceeding as
above we have
\begin{eqnarray*}
 \prod_{j=1}^n t_j^{-(a+q)}
 {\cal L}[\prod_{j=1}^n z_j^a \, f(z)]({1 \over t})
 & = &  \langle {\cal H}[{\cal K}_A(-t;x^2)](z), f(z) \rangle_L \\
 & = &  \langle {\cal K}_A(-t;z), {\cal H}[f(x^2)](z)  \rangle_L \\
 & = & {\cal L}[\prod_{j=1}^n z_j^a \, {\cal H}[f(x^2)](z)]
 \end{eqnarray*}
 Comparing this equation with (\ref{mon2}) and using the injectivity of the
 Laplace transform establishes  \linebreak (\ref{mon1}). 
 \hfill$\Box$

\subsection{Relationship to Dunkl's theory of harmonic polynomials}
In this subsection we will present the analogue of the results of 
subsection 3.5 in the type $B$ case. The required theory to do this can be 
found in ref.~\cite{dunkl89a}.

The type $B$ generalized harmonic polynomials of degree $k$ in
$x_1^2, \dots, x_n^2$ are
defined as the linearly independent solutions of the equation
$$
\Delta_B {\cal Y}_{k,l}^B = 0.
$$
Now, from \cite[Th.~1.11]{dunkl88a}, analogous to (\ref{ey}) we have
\begin{equation}\label{ey1}
E_\eta(x^2) = \sum_{m=0}^{|\eta|} r^{2m} {\cal Y}^B_{|\eta| - m,\cdot}(x^2),
\end{equation}
where $r:= (x_1^2 + \cdots + x_n^2)^{1/2}$, with
\begin{equation}
{\cal Y}^B_{|\eta| - m,\cdot}(x) =
{1 \over 4^m m! \Big ( n(a+1) + n(n-1)/\alpha +2 |\eta| - 2m \Big )_m}
\tilde{T}_{|\eta|-m}^B\Big ( {\Delta_B^m} \, E_\eta(x^2) \Big ),
\end{equation}
\begin{equation}
\tilde{T}_{k}^B = \sum_{j=0}^k{r^{2j} \over 4^j j! \Big (
-n(a+1) - n(n-1)/\alpha -2 k + 2)_j} {\Delta_B^j}.
\end{equation}

To use (\ref{ey1}) to express the non-symmetric Laguerre polynomials
in terms of the type $B$ harmonic polynomials, we note from
\cite[Prop.~3.9]{dunkl91a} that
\begin{equation}
e^{-{\Delta}_B/4} |r|^{2j} {\cal Y}^B_{k,\cdot}(x) = (-1)^j j!
L_j^{2k + n(n-1)/\alpha + n(a+1)-1}(r^2) {\cal Y}^B_{k,\cdot}(x),
\end{equation}
Thus, applying (\ref{exp.l}) to (\ref{ey1}) we obtain
\begin{equation}
E_\eta^{(L)}(x^2) = \sum_{m=0}^{|\eta|} (-1)^m m!
L_m^{(|\eta| - m +  n(n-1)/\alpha + n(a+1) -1)} (r^2)
{\cal Y}^B_{|\eta| - m,\cdot}(x),
\end{equation}
(cf.~(\ref{harmg})).

\vspace{5mm}\noindent
{\Large\bf Acknowledgments}\\[2mm]
The authors would like to thank C.~Dunkl for useful correspondence
and for communicating to us ref. \cite{dunkl96c}. The financial
support of the ARC is acknowledged.

\bibliographystyle{plain}
%\bibliography{symm}

\end{document}